\documentclass[12pt,a4paper]{article}

\usepackage[dvips]{color} 

\usepackage{amscd,amsmath,amssymb,amsfonts,times}
\usepackage{mathrsfs}
\usepackage[all]{xy}
\numberwithin{equation}{section}
    \newtheorem{thm}{Theorem}[section]
    \newtheorem{lem}[thm]{Lemma}
    
    \newtheorem{prop}[thm]{Proposition}
    \newtheorem{cor}[thm]{Corollary}

    \newtheorem{exmp}[thm]{Example}
    \newtheorem{rem}[thm]{Remark}
\newcommand{\qed}
{\mbox{}\nolinebreak$\square$\medbreak\par}
\newenvironment{pf}{\par\smallskip\noindent\emph{Proof.}}{\hfill\qed\par\smallskip}
\newenvironment{pf*}[1]{\par\smallskip\noindent\emph{#1.}}{\hfill\qed\par\smallskip}
\begin{document}
\title{Real regulator on $K_1$ of elliptic surfaces}
\author{M. Asakura}
\date\empty
\maketitle
\tableofcontents

\def\can{\omega^*}
\def\canh{\omega}
\def\cano{\mathrm{canonical}}
\def\ff{{\mathit{false}}}
\def\Coker{\mathrm{Coker}}
\def\crys{\mathrm{crys}}
\def\zar{\mathrm{zar}}
\def\dlog{d{\mathrm{log}}}
\def\dR{{\mathrm{d\hspace{-0.2pt}R}}}            
\def\et{{\mathrm{\acute{e}t}}}  
\def\Frac{{\mathrm{Frac}}}
\def\phami{\phantom{-}}
\def\id{{\mathrm{id}}}              
\def\Image{{\mathrm{Im}}}        
\def\Hom{{\mathrm{Hom}}}  
\def\ker{{\mathrm{Ker}}}          
\def\Pic{{\mathrm{Pic}}}
\def\CH{{\mathrm{CH}}}
\def\NS{{\mathrm{NS}}}
\def\NF{{\mathrm{NF}}}
\def\End{{\mathrm{End}}}
\def\pr{{\mathrm{pr}}}
\def\Proj{{\mathrm{Proj}}}
\def\ord{{\mathrm{ord}}}
\def\qis{{\mathrm{qis}}}
\def\reg{{\mathrm{reg}}}          %
\def\res{{\mathrm{res}}}          %
\def\Res{\mathrm{Res}}
\def\Spec{{\mathrm{Spec}}}     
\def\syn{{\mathrm{syn}}}
\def\cont{{\mathrm{cont}}}
\def\ind{{\mathrm{ind}}}
\def\inv{{\mathrm{inv}}}
\def\dec{{\mathrm{dec}}}
\def\Ext{{\mathrm{Ext}}}
\def\MHS{{\mathrm{MHS}}}
\def\Gr{{\mathrm{Gr}}}
\def\cand{{\mathrm{ex}}}

\def\bA{{\mathbb A}}
\def\bC{{\mathbb C}}
\def\C{{\mathbb C}}
\def\G{{\mathbb G}}
\def\bE{{\mathbb E}}
\def\bF{{\mathbb F}}
\def\F{{\mathbb F}}
\def\bH{{\mathbb H}}
\def\bJ{{\mathbb J}}
\def\bN{{\mathbb N}}
\def\bP{{\mathbb P}}
\def\P{{\mathbb P}}
\def\bQ{{\mathbb Q}}
\def\Q{{\mathbb Q}}
\def\bR{{\mathbb R}}
\def\R{{\mathbb R}}
\def\bZ{{\mathbb Z}}
\def\Z{{\mathbb Z}}
\def\cA{{\mathscr A}}
\def\cD{{\mathscr D}}
\def\cM{{\mathscr M}}
\def\cL{{\mathscr L}}
\def\cE{{\mathscr E}}
\def\cO{{\mathscr O}}
\def\O{{\mathscr O}}
\def\cR{{\mathscr R}}
\def\cS{{\mathscr S}}
\def\cX{{\mathscr X}}
\def\cH{{\mathscr H}}
\def\ccH{{\mathscr H}_C}
\def\PF{{\mathit{PF}}}
\def\Div{{\mathrm{Div}}}
\def\codim{{\mathrm{codim}}}
%
\def\ep{\epsilon}
\def\vG{\varGamma}
\def\vg{\varGamma}

%
%
%
%
\def\lra{\longrightarrow}
\def\lla{\longleftarrow}
\def\Lra{\Longrightarrow}
\def\hra{\hookrightarrow}
\def\lmt{\longmapsto}
\def\ot{\otimes}
\def\op{\oplus}
\def\wt#1{\widetilde{#1}}
\def\wh#1{\widehat{#1}}
\def\spt{\sptilde}
\def\ol#1{\overline{#1}}
\def\ul#1{\underline{#1}}
\def\us#1#2{\underset{#1}{#2}}
\def\os#1#2{\overset{#1}{#2}}
\def\lim#1{\us{#1}{\varinjlim}}
\def\plim#1{\us{#1}{\varprojlim}}

\section{Introduction}
Let $X$ be a projective nonsingular
variety over the complex number field $\C$.
Let $H^i_\cM(X,\Z(j))$ denotes the motivic cohomology group. 
It is known that 
$H^i_\cM(X,\Q(j))$ is isomorphic to Quillen's $K$-group $K_{2j-i}(X)^{(j)}$.
By the theory of higher Chern classes,
we have the {\it Beilinson regulator map} (higher Chern class map)
\[
\reg_{i,j}:H^i_\cM(X,\Z(j))\lra H_\cD^i(X,\Z(j))
\]
to the Deligne-Beilinson cohomology group
(\cite{schneider}).
The purpose of this paper is to give a certain method for computations
of the regulator map for $(i,j)=(3,2)$ (namely $K_1$)
and $X$ an elliptic surface.

\medskip

The cup-product pairing gives rise to a map 
$\C^\times\ot \Pic(X)\cong \C^\times\ot H^2_\cM(X,\Z(1))
\to H^3_\cM(X,\Z(2))$. Its image is called the {\it decomposable} part,
and the cokernel is called the {\it indecomposable} part.
The decomposable part does not affect serious difficulty,
while the indecomposable part plays the central role in the study of
$H^3_\cM(X,\Z(2))$.
According to \cite{lewisJAG}, we call an element $\xi\in H^3_\cM(X,\Z(2))$
{\it regulator indecomposable} if
$\reg_{3,2}(\xi)$ does not lie in the image of $\C^\times\ot\NS(X)$.
Obviously regulator indecomposable elements are indecomposable.
The converse is also true if the Beilinson-Hodge conjecture for $K_2$ is true.
Lewis and Gordon constructed regulator indecomposable elements in case
$X$ is a product of `general' elliptic curves (\cite{lewisJAG} Theorem 1).
There are a lot of other related works, though I don't catch up all of them.
On the other hand, in case that $X$ is defined over a number field,
the question is more difficult, and as far as I know
there are only a few of such examples (e.g. \cite{R} \S 12).

The real regulator map $\reg_{3,2}$ is usually written
in terms of differential $(1,1)$-forms.
Then one of the technical difficulties appears from the fact that it is not easy
to describe analytic differential forms explicitly.
The key idea in this paper is to use
certain ``algebraic" 2-forms instead of analytic forms. 
This makes it easier to describe and compute the real regulator.

\medskip

This paper is organized as follows.
\S \ref{reg-sect} is a quick review of $H^3_\cM(X,\Q(2))$ and Beilinson regulator.
In \S \ref{deRh-sect}, we provide notations and some elementary results on 
de Rham cohomology and the Hodge filtration. Especially we introduce ``good algebraic 2-forms"
which plays a key role in our computations (\S \ref{good-sect}).
In \S \ref{exp-sect}, we give a method of computations of real regulator on
$K_1$ of elliptic surfaces.
In \S \ref{Example-sect} we give an example. In particular
we construct regulator indecomposable elements for an elliptic surface
defined over $\Q$ with arbitrary large $p_g$ (Cor.\ref{Example-cor}).
\S \ref{Appendix-sect} is an appendix providing proofs 
of some explicit formulas on Gauss-Manin connection.

\bigskip

\noindent{\bf Acknowledgment.}
A rough idea was inspired
during my visit to the University of Alberta in September 2012,
especially when I discussed the paper \cite{lewis}
 with Professor James Lewis.
I would like to express special thanks to him.
I'd also like to thank the university members
for their hospitality.

\section{Real Regulator map on $H^3_\cM(X,\Q(2))$}\label{reg-sect}
For a regular and integral scheme $X$,
let $Z_i(X)=Z^{\dim X-i}(X)$ be the free abelian group of irreducible subvarieties
of Krull dimension $i$.
For an integral scheme $X$, we denote by $\eta_X$ the field of rational functions on $X$.
For schemes $X$ and $T$ over a base scheme $S$, we set $X(T)=\mathrm{Mor}_S(T,X)$.
and say $x\in X(T)$ a $T$-valued point of $X$.
If $T=\Spec R$, then we also write $X(R)=X(\Spec R)$.

\subsection{$H^3_\cM(X,\Q(2))$ and (in)decomposable parts}
Let $X$ be a smooth variety over a field $K$.
Let $D\subset X$ be an irreducible divisor, and $\wt{D}\to D$
the normalization.
Let $j:\wt{D}\to D\hra X$ be the composition.
Then 
we define $\Div_D(f):=j_*\Div_{\wt{D}}(f)\in Z^2(X)$ the push-forward of the Weil divisor
on $\wt{D}$ by $j$.
Let
\[
\partial_1:\bigoplus_{\codim D=1} \eta_D^\times\lra Z^2(X),\quad
[f,D]\longmapsto \Div_D(f)
\]
be a homomorphism where we write
\[
[f,D]:=(\cdots,1,f,1,\cdots)\in \bigoplus_{\codim D=1} \eta_D^\times,\quad
\mbox{($f$ is placed in the $D$-component).}
\] 
Let \[\partial_2:K^M_2(\eta_X)\lra \bigoplus_{\codim D=1} \eta_D^\times\]
\[
\partial_2\{f,g\}=\sum_{\codim D=1}
\left[(-1)^{\ord_D(f)\ord_D(g)}\frac{f^{\ord_D(g)}}{g^{\ord_D(f)}}|_D,D\right]
\]
be the tame symbol.
Then it is well-known that there is the canonical isomorphism
\begin{equation}\label{mot-0}
H^3_\cM(X,\Q(2))\cong \left(\frac{\ker(\bigoplus \eta_D^\times
\os{\partial_1}{\lra} Z^2(X))}
{\Image(K^M_2(\eta_X)\os{\partial_2}{\lra} \bigoplus \eta_D^\times}\right)\ot\Q.
\end{equation}
In this paper we always identify the motivic cohomology group $H^3_\cM(X,\Q(2))$
with the group in the right hand side of \eqref{mot-0}.

Let $L/K$ be a finite extension. Write $X_L:=X\times_KL$.
Then there is the obvious map
\[
L^\times\ot Z^1(X_L)\lra H^3_\cM(X_L,\Q(2)),\quad \lambda\ot D\longmapsto[\lambda,D].
\]
Let $N_{L/K}:H^3_\cM(X_L,\Q(2))\to H^3_\cM(X,\Q(2))$ be the norm map on motivic 
cohomology.
Then we put
\[
H^3_\cM(X,\Q(2))_\dec:=\sum_{[L:K]<\infty}
N_{L/K}(\Image(L^\times\ot Z^1(X_L)\to H^3_\cM(X_L,\Q(2))))
\]
and call it the {\it decomposable part}.
We put
\[
H^3_\cM(X,\Q(2))_\ind:=
H^3_\cM(X,\Q(2))/H^3_\cM(X,\Q(2))_\dec
\]
and call it the {\it indecomposable part}.
The indecomposable part plays the central role in the study of $H^3_\cM(X,\Q(2))$.

\subsection{Beilinson regulator on indecomposable parts}\label{bregin-sect}
For a smooth projective variety $X$ over $\C$,
we denote by $H^\bullet_B(X,\Q)=H^\bullet_B(X(\C),\Q)$ 
(resp. $H_\bullet(X,\Q)$) the Betti cohomology (resp. Betti homology).
$H^\bullet_\dR(X)=H^\bullet_\dR(X/\C)$ denotes the de Rham cohomology.

\medskip

By the theory of universal Chern class, there is the {\it Beilinson regulator map}
\begin{align}\label{beireg-0}
\reg=\reg_\Q:H^3_\cM(X,\Q(2))\lra H^3_\cD(X,\Q(2))&\cong\Ext_\MHS^1(\Q,H^2(X,\Q(2)))\\
&=\frac{H^2_B(X,\C)}{F^2+H^2_B(X,\Q(2))}
\end{align}
to the Deligne-Beilinson cohomology group, which is isomorphic to
the Yoneda extension group of mixed Hodge structures
where $H^2(X,\Q(2))=(H^2_B(X,\Q(2)),F^\bullet H^2_\dR(X))$ denotes the Hodge structure
(of weight $-2$).
Put
\[
H^2_B(X)_\ind:=H^2_B(X,\Q(1))/\NS(X)\ot\Q,\quad
H^2_\dR(X)_\ind:=H^2_\dR(X)/\NS(X)\ot\C,
\]
\[
H^2(X)_\ind:=(H^2_B(X)_\ind,F^\bullet H^2_\dR(X)_\ind)
\mbox{ (= a Hodge structure of weight 0).}
\]
Then \eqref{beireg-0} yields a commutative diagram
\begin{equation}\label{bregin-2}
\xymatrix{
0\ar[d]&0\ar[d]\\
H^3_\cM(X,\Q(2))_\dec\ar[d]\ar[r]&\Ext_\MHS^1(\Q,\NS(X)\ot\Q(1)))\ar[d]\\
H^3_\cM(X,\Q(2))\ar[d]\ar[r]^{\reg\qquad}&\Ext^1_\MHS(\Q,H^2(X,\Q(2)))\ar[d]\\
H^3_\cM(X,\Q(2))_\ind\ar[d]\ar[r]^{\ol{\reg}\quad\qquad}&\Ext_\MHS^1(\Q,H^2(X)_\ind\ot\Q(1))\ar[d]\\
0&0
}
\end{equation}
The top arrow is simply written by ``$\log$", namely the composition
\[
\C^\times\ot\Pic(X)\to H^3_\cM(X,\Q(2))_\dec\lra 
\Ext_\MHS^1(\Q,\NS(X)\ot\Q(1)))\cong \C/\Q(1)\ot\NS(X)
\]
is given by $\lambda\ot Z\mapsto \log(\lambda)\ot Z$.
The bottom arrow $\ol{\reg}$ plays an important role.
Let us describe it in terms of extension of mixed Hodge structures.
Let $n=\dim  X$.
Let $\xi=\sum[f_i,D_i]\in \bigoplus\eta_{D_i}^\times$ such that $\partial_1(\xi)=0$.
Let $\reg'$ be the composition
\begin{align*}
H^3_\cM(X,\Q(2))&\os{\reg}{\lra}\Ext^1_\MHS(\Q,H^2(X,\Q(2)))\\
&\lra \Ext^1_\MHS(\Q,H^2(X,\Q(2))/\langle D_i\rangle)\\
&\os{\cong}{\lra}
\Ext^1_\MHS(\Q,H_{2n-2}(X,\Q(2-n))/H_{2n-2}(D,\Q(2-n)))
\end{align*}
where $\langle D_i\rangle$ denotes the subgroup generated by the cycle classes of
$D_i$, and the last isomorphism is the Poincare duality.
Let $j:\wt{D}_i\to D_i$ be the normalization.
Let $\wt{Z}_i\subset \wt{D}_i$ be the support of $\Div_{\wt{D}_i}(f_i)$. 
Put
\[
\wt{Z}:=\coprod_i \wt{Z}_i\subset \wt{D}:=\coprod_i \wt{D}_i,
\quad
Z:=\bigcup_i j(\wt{Z}_i)\subset D:=\bigcup_i D_i.
\]
Consider a commutative diagram
\[
\xymatrix{
&&H^1(\wt{D}-\wt{Z},\Z(1))\ar[d]^a\\
0\ar[r]&H_{2n-3}(\wt{D},\Z(2-n))\ar[d]^{j_*}\ar[r]&
H_{2n-3}(\wt{D},\wt{Z};\Z(2-n))\ar[d]^{j_*}\ar[r]^{\delta_1}
&H_{2n-4}(\wt{Z},\Z(2-n))\ar[d]^{j_*}\\
0\ar[r]&H_{2n-3}(D,\Z(2-n))\ar[r]^b&H_{2n-3}(D,Z;\Z(2-n))\ar[r]^{\delta_2}
&H_{2n-4}(Z,\Z(2-n))
}
\]
with exact rows.
Let
\[
\nu:=\left(\frac{df_i}{f_i}\right)\in H^1(\wt{D}-\wt{Z},\Z(1)).
\]
Since $\partial_1(\xi)=0$, one has $j_*\delta_1a(\nu)=0$.
Therefore $\nu$ defines
$\nu_\xi\in H_{2n-3}(D,\Z(2-n))$ such that $b(\nu_\xi)=j_*a(\nu)$.
Note that $\nu_\xi$ belongs to the Hodge (0,0)-part because so does $\nu$.
By the exact sequence
\[
\cdots\lra H_{2n-2}(X,D;\Q(2-n))\os{\partial}{\lra} H_{2n-3}(D,\Q(2-n))\os{\delta}{\lra} 
H_{2n-3}(X,\Q(2-n))\lra \cdots
\]
we have an exact sequence
\begin{equation}\label{beireg-ex}
0\to H_{2n-2}(X,\Q(2-n))/H_{2n-2}(D,\Q(2-n))\to H_{2n-2}(X,D;\Q(2-n))
\os{\partial}{\to}
\ker(\delta)\to0
\end{equation}
of mixed Hodge structures.
Since the weight of $H_{2n-3}(X,\Q(2-n))$ is $-1$,  
the Hodge $(0,0)$-part of $H_{2n-3}(D,\Q(2-n))$ is contained in the kernel of $\delta$.
In particular 
we have 
an exact sequence 
\begin{equation}\label{beireg-1}
0\lra H_{2n-2}(X,\Q(2-n))/H_{2n-2}(D,\Q(2-n))\to H_\xi(X,D)\lra
\Q\lra 0
\end{equation}
by taking the pull-back of \eqref{beireg-ex} via $\Q\to \ker(\delta)$, $1\mapsto \nu_\xi$.
Then the following is well-known to specialists, proven by using the Riemann-Roch
theorem without denominators
(\cite{gillet}, see also \cite{sato} Thm. 11.2).
\begin{thm}\label{bregin-thm}
$\reg'(\xi)$ corresponds to \eqref{beireg-1} up to sign.
In other words, letting
\[
\rho:\Q\lra \Ext^1_{\MHS}(\Q,H_{2n-2}(X,\Q(2-n))/H_{2n-2}(D,\Q(2-n)))
\]
be the connecting homomorphism arising from \eqref{beireg-1}, one has $\reg'(\xi)=\pm\rho(1)$.
\end{thm}
For the later use, we write down $\rho(1)$ explicitly.
Write \[
M:=H_{2n-2}(X,\Q(2-n))/H_{2n-2}(D,\Q(2-n)),\quad
H^{2n-2}_\dR(X)':=\ker[ H^{2n-2}_\dR(X)\lra H^{2n-2}_\dR(D)].\]
Then the natural isomorphism
$M\ot_\Q\C\cong
\Hom(H^{2n-2}_\dR(X)',\C)
$
induces
\begin{equation}\label{beireg-2}
\Ext^1_{\MHS}(\Q,M)
\cong
\Coker[
H_{2n-2}(X,\Q(2-n))\os{\Phi}{\to} \Hom(F^{n-1}H^{2n-2}_\dR(X)',\C)
]
\end{equation}
where
\[
\Phi(\Delta)=\left[
\omega\longmapsto \int_\Delta\omega
\right],\quad \omega\in F^{n-1}H^{2n-2}_\dR(X)'.
\]
Taking the dual of the map $\Q\to H_{2n-3}(D,\Q(2-n))$,
$1\mapsto \nu_\xi$, one has $H^{2n-3}_\dR(D)\to \C$
and this induces
\[
0\lra \C\lra H^{2n-2}_\dR(X,D)'\lra H^{2n-2}_\dR(X)'\lra 0,
\]
which is isomorphic to the dual of \eqref{beireg-1}.
Let $\omega_{X,D}\in F^{n-1}H^{2n-2}_\dR(X,D)'$ denotes the element corresponding to
$\omega\in F^{n-1}H^{2n-2}_\dR(X)'$ via the isomorphisms
$F^{n-1}H^{2n-2}_\dR(X,D)'\os{\cong}{\to} F^{n-1}H^{2n-2}_\dR(X)'$.
Let $\Gamma\in H_{2n-2}(X,D;\Q(2-n))$ be an arbitrary element 
such that $\partial(\Gamma)=\nu_\xi$.
Then we have
\begin{equation}\label{beireg-3}
\rho(1)=\left[
\omega\longmapsto \int_\Gamma\omega_{X,D}
\right]
\end{equation}
under the isomorphism \eqref{beireg-2}.

\medskip

The {\it real regulator map} is the composition of $\reg_\Q$ and the canonical map
\[
\Ext^1_\MHS(\Q,H^2(X,\Q(2)))\lra
\Ext^1_{\R\mbox{-}\MHS}(\R,H^2(X,\R(2)))
\]
to the extension group of real mixed Hodge structures, 
which we denote by $\reg_\R$:
\begin{equation}\label{bregin-3}
\reg_\R:H^3_\cM(X,\Q(2))\lra \Ext_{\R\mbox{-}\MHS}^1(\R,H^2(X,\R(2)))
\cong H^2_B(X,\R(1))\cap H^{1,1}.
\end{equation}
This also induces 
\begin{equation}
\ol{\reg}_\R:H^3_\cM(X,\Q(2))_\ind\to \Ext_{\R\mbox{-}\MHS}^1(\R,H^2(X)_\ind\ot\R(1))
\cong (H^2_B(X)_\ind\ot\R)\cap H^{1,1}
\end{equation}
on the indecomposable part.

\subsection{$\Q$-structure on determinant of $H^3_\cD(X/\R,\R(2))$}
Suppose that $X$ is a projective smooth variety over $\Q$.
Write $X_\C:=X\times_\Q\C$.
The {\it infinite Frobenius} map $F_\infty$ is defined to be the
anti-holomorphic map on $X(\C)=\mathrm{Mor}_\Q(\Spec\C,X)$ 
induced from the complex conjugation on $\Spec\C$.
For a subring $A\subset \R$,
the infinite Frobenius map 
acts on the Deligne-Beilinson complex $A_X(j)_\cD$ in a canonical way,
so that we have the involution on $H^\bullet_\cD(X_\C,A(j))$, which we denote by the
same notation $F_\infty$.
We define
\[
H^\bullet_\cD(X/\R,A(j)):=H^\bullet_\cD(X_\C,A(j))^{F_\infty=1}
\]
the fixed part by $F_\infty$. We call it the {\it real Deligne-Beilinson cohomology}.
Since the action of $F_\infty$ is compatible via the Beilinson regulator map,
we have
\begin{align}
\reg_\R:H^3_\cM(X,\Q(2))\lra H_{\cD}:=&
\Ext_{\R\mbox{-}\MHS}^1(\R,H^2(X_\C,\R(2)))^{F_\infty=1}\\
\cong& \frac{H^2_B(X_\C,\R(1))^{F_\infty=1}}{F^2H^2_\dR(X/\R)},\label{qstr-3}
\end{align}
and 
\begin{align}
\ol{\reg}_\R:H^3_\cM(X,\Q(2))\lra H_{\cD,\ind}:=&
\Ext_{\R\mbox{-}\MHS}^1(\R,H^2(X_\C)_\ind\ot\R(1))^{F_\infty=1}\\
\cong& \frac{[H^2_B(X_\C)_\ind\ot\R]^{F_\infty=1}}{F^2H^2_\dR(X/\R)}.\label{qstr-4}
\end{align}
There are the canonical $\Q$-structures $e_\Q$ and $e_{\ind,\Q}$ 
on the determinant vector spaces $\det H_\cD$ and
$\det H_{\cD,\ind}$:
\[
\R\cdot e_\Q=\det H_\cD,\quad
\R\cdot e_{\ind,\Q}=\det H_{\cD,\ind}.
\]
Here we recall the definition.
The isomorphisms \eqref{qstr-3} and \eqref{qstr-4} induce 
\begin{equation}\label{qstr-1}
\det H_\cD
\cong 
\det [H^2_B(X_\C,\R(1))^{F_\infty=1}]\ot[\det F^2H^2_\dR(X/\R)]^{-1},
\end{equation}
and 
\begin{equation}\label{qstr-2}
\det H_{\cD,\ind}
\cong \det[(H^2_B(X_\C)_\ind\ot\R)^{F_\infty=1}]
\ot[\det F^2H^2_\dR(X/\R)]^{-1}.
\end{equation}
The right hand sides of \eqref{qstr-1} and \eqref{qstr-2} have the $\Q$-structures
induced from the $\Q$-structures
\[
H^2_B(X_\C,\Q(1))^{F_\infty=1},\quad
H^2_B(X_\C)_\ind^{F_\infty=1}, \quad F^2H^2_\dR(X/\Q).
\]
The $\Q$-structures $e_\Q$ and $e_{\ind,\Q}$ are defined to be the corresponding
one:
\begin{equation}\label{qstr-9}
\Q\cdot e_\Q\cong 
\det [H^2_B(X_\C,\Q(1))^{F_\infty=1}]\ot[\det F^2H^2_\dR(X/\Q)]^{-1},
\end{equation}
\begin{equation}\label{qstr-10}
\Q\cdot e_{\Q,\ind}\cong 
\det[H^2_B(X_\C)_\ind^{F_\infty=1}]
\ot[\det F^2H^2_\dR(X/\Q)]^{-1}.
\end{equation}
\subsection{$e_\Q^\ff$ and $e_{\ind,\Q}^\ff$}\label{false-sect}
We introduce other $\Q$-structures $e_\Q^\ff$ and $e_{\ind,\Q}^\ff$ on $\det H_\cD$ and
$\det H_{\cD,\ind}$. 
For simplicity, we assume $\dim X=2$.
Put
\[
H_2(X_\C,\Q)_\ind:=H_2(X_\C,\Q)/(\NS(X_\C)\ot\Q(1))\cong 
H^2_B(X_\C)_\ind\ot\Q(1),
\]
\[
H^2_\dR(X/\Q)_\ind:=\mathrm{Coim}(H^2_\dR(X/\Q)\lra H^2_\dR(X/\C)/(\NS(X_\C)\ot\C)).
\]
Note that $H^2_\dR(X/\Q)_\ind\ot\C\os{\cong}{\to}
H^2_\dR(X/\C)/(\NS(X_\C)\ot\C)$.
There are exact sequences
\begin{equation}\label{qstr-5}
0\lra H_2(X_\C,\R)^{F_\infty=1}\lra \Hom(F^1H^2_\dR(X/\Q),\R)\lra H_\cD\lra 0
\end{equation}
\begin{equation}\label{qstr-6}
0\lra H_2(X_\C,\R)_\ind^{F_\infty=1}\lra \Hom(F^1H^2_\dR(X/\Q)_\ind,\R)\lra 
H_{\cD,\ind}\lra 0
\end{equation}
under the canonical isomorphisms
\begin{equation}\label{qstr-7}
H^2_B(X_\C,\C)\cong H^2_\dR(X/\C),\quad H^2_B(X_\C,\Q(2))\cong H_2(X_\C,\Q).
\end{equation}
Then the $\Q$-structures
\[
H_2(X_\C,\Q)^{F_\infty=1},\quad H_2(X_\C,\R)_\ind^{F_\infty=1},\quad
H^2_\dR(X/\Q),\quad
H^2_\dR(X/\Q)_\ind,
\]
induce $e^{\ff}_\Q$ and $e^\ff_{\ind,\Q}$:
\begin{equation}\label{qstr-11}
\Q\cdot e_\Q^\ff\cong 
[\det H_2(X_\C,\Q)^{F_\infty=1}]^{-1}\ot[\det F^1H^2_\dR(X/\Q)]^{-1},
\end{equation}
\begin{equation}\label{qstr-12}
\Q\cdot e_{\ind,\Q}^\ff\cong 
[\det H_2(X_\C,\Q)_\ind^{F_\infty=1}]^{-1}\ot[\det F^1H^2_\dR(X/\Q)_\ind]^{-1}.
\end{equation}

\begin{prop}\label{qstr-8}
Put\[
r:=\dim H_2(X_\C,\Q)^{F_\infty=1}=\dim H^2_B(X_\C,\Q(1))^{F_\infty=-1},
\] 
\[
s:=\dim H_2(X_\C,\Q)_\ind^{F_\infty=1}=\dim H^2_B(X_\C)_\ind^{F_\infty=-1}
=r-\dim \NS(X_\C)^{F_\infty=-1}.
\]
Write \[
H_B:=H^2_B(X_\C,\Q(1)), \quad H_{B,\ind}:=H^2_B(X_\C)_\ind,\] 
\[
F^\bullet H_\dR:=F^\bullet H^2_\dR(X/\Q),\quad
F^\bullet H_{\dR,\ind}:=F^\bullet H^2_\dR(X/\Q)_\ind\] simply.
Then
\[
\Q\cdot e^{\ff}_\Q=\Q\cdot e_\Q\ot\Q(-r)\ot \det H_\dR\ot
[\det H_B]^{-1},
\]
\[
\Q\cdot e^{\ff}_{\ind,\Q}=\Q\cdot e_{\ind,\Q}\ot\Q(-s)\ot \det H_{\dR,\ind}\ot
[\det H_{B,\ind}]^{-1},
\]
where we mean
\[
\det H_\dR\ot
[\det H_B]^{-1}\subset \det H_\dR^2(X/\C)\ot
[\det H_B^2(X_\C,\C)]^{-1}\os{\eqref{qstr-7}}{\cong}\C,\mbox{ etc.}\]
\end{prop}
\begin{pf}
By the Poincare duality,
\[
\det F^1H_\dR=[\det H_\dR]^{-1}\ot \det F^2H_\dR,\quad
\det F^1H_{\dR,\ind}=[\det H_{\dR,\ind}]^{-1}\ot \det F^2H_\dR.
\]
Moreover
\[
\det [H_2(X_\C,\Q)^{F_\infty=1}]
=\det [H^2_B(X_\C,\Q(2))^{F_\infty=1}]
=\det [H_B^{F_\infty=-1}\ot\Q(1)]
=\Q(r)\ot \det H_B^{F_\infty=-1}
\]
and
\[
\det [H_2(X_\C,\Q)_\ind]^{F_\infty=1}
=\det[ H_{B,\ind}^{F_\infty=-1}\ot\Q(1)]
=\Q(s)\ot\det H_{B,\ind}^{F_\infty=-1}.
\]
Therefore we have
\begin{align*}
\Q\cdot e^{\ff}_\Q\ot e_\Q^{-1}
&=
\Q(-r)\ot[\det H_B^{F_\infty=-1}]^{-1}
\ot [\det H_B^{F_\infty=1}]^{-1}\ot[\det H_\dR]\\
&=
\Q(-r)\ot[\det H_B]^{-1}
\ot[\det H_\dR]\\
\end{align*}
by \eqref{qstr-9} and \eqref{qstr-11},
and
\begin{align*}
\Q\cdot e^{\ff}_{\ind,\Q}\ot e_{\ind,\Q}^{-1}
&=
\Q(-s)\ot[\det H_{B,\ind}^{F_\infty=-1}]^{-1}
\ot [\det H_{B,\ind}^{F_\infty=1}]^{-1}\ot[\det H_{\dR,\ind}]\\
&=
\Q(-s)\ot[\det H_{B,\ind}]^{-1}
\ot[\det H_{\dR,\ind}]\\
\end{align*}
by \eqref{qstr-10} and \eqref{qstr-12}. This completes the proof.
\end{pf}
\begin{rem}
The Poincare duality implies
\[
(\det H_B)^{\ot 2}\cong 
 H^4_B(X_\C,\Q(2))^{\ot m}\cong
H^4_\dR(X/\Q)^{\ot m}\cong
(\det H_\dR)^{\ot 2},
\]
and
\[
(\det H_{B,\ind})^{\ot 2}\cong 
 H^4_B(X_\C,\Q(2))^{\ot m'}\cong
H^4_\dR(X/\Q)^{\ot m'}\cong
(\det H_{\dR,\ind})^{\ot 2}.
\]
Therefore $(\det H_\dR\ot[\det H_B]^{-1})$ and
$(\det H_{\dR,\ind}\ot[\det H_{B,\ind}]^{-1})$ are contained in $\sqrt{\Q^\times}$
(possibly rational numbers).
\end{rem}


\section{Elliptic surface
and good algebraic 2-forms}\label{deRh-sect}
\subsection{Notations}\label{deRh-sect-1}
Let $K$ be a field of characteristic 0.
Let $f:X\to C$ be an elliptic surface with a section $e:C\to X$.
This means that $X$ (resp. $C$) is a projective smooth surface (resp. curve) over $K$,
and the generic fiber of $f$ is an elliptic curve.
Hereafter we assume that the $j$-invariant of $f$ is not constant, namely,
$f$ is not isotrivial.

\medskip

Throughout \S \ref{deRh-sect} and \S \ref{exp-sect}
we use the following notations.
\begin{itemize}
\item
$D\subset X$ is the sum of the multiplicative fibers. Put $T_m=f(D)\subset C$.
Note $T_m\ne\emptyset$ by the assumption.
\item
$E\subset X$ is the sum of the additive fibers. Put $T_a=f(E)\subset C$.
\item
$S=C-(T_m+T_a)$ and $U=f^{-1}(S)=X-(E+D)$.
\item
$\ol{S}=C-T_a$ and $\ol{U}=f^{-1}(\ol{S})=X-E$.
\item
Let $F\subset S$ be the support of the cokernel of the $\O_S$-linear map 
\begin{equation}\label{PF-rem-eq}
\ol{\nabla}:f_*\Omega^1_{U/S}\lra \Omega^1_{S}\ot R^1f_*\O_{U}
\end{equation}
induced from the Gauss-Manin connection.
(By Cor. \ref{PF-rem}, this is a set of finite closed points.)
Hence $\ol{\nabla}$ is an isomorphism outside $F$.
\item
Put $S^o:=S-F=C-(T_m+T_a+F)$ and $U^o:=f^{-1}(S^o)$.
\item
$\ol{S^o}:=S^o+T_m=C-(T_a+F)$
and 
$\ol{U^o}:=f^{-1}(\ol{S})$.
\item
Write $X_{\ol{K}}=X\times_K\ol{K}$.
Let $\NF(X_{\ol{K}})\subset \NS(X_{\ol{K}})$ denotes the subgroup of the Neron-Severi group
generated by $e(C)$ and irreducible components of $D_{\ol{K}}+E_{\ol{K}}$.
\item
$\NF_\dR(X)=H^2_\dR(X/K)\cap (\NF(X_{\ol{K}})\ot_\Z\ol{K})\subset H^2_\dR(X_{\ol{K}}/\ol{K})$.
\end{itemize}
\[
\xymatrix{
U\ar[r]\ar[d]\ar@{}[rd]|{\square}
&\ol{U}\ar[r]\ar[d]\ar@{}[rd]|{\square}
&X\ar[d]^f&\ol{U^o}\ar@{}[ld]|{\square}\ar[d]\ar[l]
&U^o\ar@{}[ld]|{\square}\ar[d]\ar[l]&S\ar[r]&\ol{S}\\
S\ar[r]&\ol{S}\ar[r]&C&\ol{S^o}\ar[l]&S^o\ar[l]&S^o=\ol{S^o}\cap S\ar[u]\ar[r]
&\ol{S^o}\ar[u]
}
\]

\begin{rem}\label{nf-rem-dege}
The intersection pairing $\NF(X_{\ol{K}})\ot \NF(X_{\ol{K}})\to \Q$ is non-degenerate.
This is proven on a case-by-case analysis by using the classification of degenerations
(see \cite{Si} IV, Thm.8.2 for the classification).
\end{rem}
\begin{rem}\label{nf-rem}
$\NF_\dR(X)\ot_K\ol{K}=\NF(X_{\ol{K}})\ot_\Z\ol{K}$ in $H^2_\dR(X_{\ol{K}}/\ol{K})$.
This is proven by using \cite{AEC} II Lemma 5.8.1.
\end{rem}
\begin{rem}
By Cor. \ref{PF-rem}, $F$ is described in the following way.
Around a neighborhood of $s\in S$, $f$ is written by a Weierstrass form $y^2=4x^3-g_2x-g_3$
$(g_2,g_3\in\O_{S,s}$, $\Delta=g_2^3-27g_3^2\in\O_{S,s}^\times$).
Let $j=1728g_2^3/(g_2^3-27g_3^2)$ be the $j$-invariant. Then $s\in F$ if and only if
\[
\frac{g_2}{g_3}\frac{dj}{j}\in \O_{S,s}\ot\Omega^1_S\cong\O_{S,s}
\quad (\Omega^1_S:=\Omega^1_{S/K})
\]
is a free $\O_{S,s}$-basis. 
\end{rem}

\begin{prop}\label{VDR}
Let $Q\subset C$ be a non-empty open set, and $V:=f^{-1}(Q)$.
We put
\[
H^2_\dR(V)_0:=\ker[H^2_\dR(V)\lra \prod_{s\in Q}H^2_\dR(f^{-1}(s))\times H^2_\dR(V\cap
e(C))].
\]
When $V=X$, we also write $\NF_\dR(X)^\perp=H^2_\dR(X)_0(=$ the
orthogonal complements of 
$\NF_\dR(X)$ in $H^2_\dR(X)$
with respect to the cup-product pairing$)$. 
Then the following hold.
\begin{enumerate}
\item[{\rm (1)}]
If $V\ne X$, then $H^2_\dR(V)_0=\Image[\vg(V,\Omega^2_V)\lra H^2_\dR(V)]$.
\item[{\rm (2)}] Let $Q_1\supset Q_2$ and $V_i=f^{-1}(Q_i)$. 
Then there is an exact sequence
\[
0\lra H^2_\dR(V_1)_0\lra H^2_\dR(V_2)_0\lra \bigoplus_{s\in Q_1-Q_2}H_{1,\dR}(f^{-1}(s)).
\]
\item[{\rm (3)}]  
$\NF_\dR(X)^\perp\os{\cong}{\to} H^2_\dR(\ol{U})_0$.
\end{enumerate}
\end{prop}
\begin{pf}
Note that $\NF_\dR(X)^\perp\ot_K\ol{K}=(\NF(X_{\ol{K}})\ot_\Z\ol{K})^\perp$ 
by Rem. \ref{nf-rem}.
Therefore we may assume $K=\ol{K}$ throughout the proof pf Prop. \ref{VDR}.

\medskip

We consider a spectral sequence
\[
E^{pq}_1=H^q(V,\Omega^p_{V})\Longrightarrow
H^{p+q}_\dR(V).
\]
Since $Q$ is affine by the assumption, $E_1^{pq}=H^q(V,\Omega^p_{V})
=\vg(Q,R^qf_*\Omega^p_{V})=0$ unless $p\leq 2$ and $q\leq 1$, so that
we have
\[
E^{20}_3=E^{20}_\infty=\Image\vg(V,\Omega^2_{V}),\quad E^{11}_2=E^{11}_\infty,
\quad E_2^{02}=0,
\]
\[
0\lra \Image\vg(V,\Omega^2_{V})\lra H^2_\dR(V)\lra E^{11}_\infty\lra0. 
\]
\begin{lem}\label{VDR-lem1}
$E^{11}_2=E^{11}_\infty$ is generated by the image of the cycle classes 
of $e(C)$ and irreducible components of each fiber $f^{-1}(s)$ as $K$-module
(note we assumed $K=\ol{K}$ throughout the proof).
\end{lem}
\begin{pf}
Let $Q^o:=Q\cap S^o$ and $j:V^o:=f^{-1}(Q^o)\hra V$ be the open immersion.
Consider a commutative diagram
\[
\xymatrix{
\vg(V,j_*\Omega^1_{V^o}/\Omega^1_{V})\ar[r]^{\quad \delta}&
H^1(V,\Omega^1_{V})\ar[r]^{j^*}\ar[d]^{d}&
H^1(V^o,\Omega^1_{V^o})\ar[d]^{d}&\mbox{(exact)}\\&
H^1(V,\Omega^2_{V})\ar[r]&
H^1(V^o,\Omega^2_{V^o})
}
\]
Let $x\in \ker d$. 
Then $j^*(x)\in \ker d$.
The kernel of $d$ on $H^1(V^o,\Omega^1_{V^o})$ is generated by the cycle class 
$[e(C)]$ as $K$-module.
Indeed, since $\ol{\nabla}$ \eqref{PF-rem-eq} is an isomorphism on $V^o$, 
one has $R^1f_*\Omega^1_{V^o}\os{\cong}{\lra}R^1f_*\Omega^1_{V^o/Q^o}$, 
and this is generated by the cycle class of $e(C)$ as $\O_{Q^o}$-module.
Then one can identify the map $d$ with $d\ot\id:\O(Q^o)\ot[e(C)]\to
\vg(Q^o,R^1f_*\Omega^2_{V^o})\cong
\vg(Q^o,\Omega^1_{Q^o}\ot R^1f_*\Omega^1_{V^o/Q^o})
=\vg(Q^o,\Omega^1_{Q^o})\ot [e(C)]$. Since the characteristic of $K$ is zero, the kernel of it
is one-dimensional over $K$. This means $\ker~d$ is generated by the cycle class $[e(C)]$.
Thus $x':=x-c[e(C)]$ for some $c\in K$
is contained in $\ker(j^*)=\Image~\delta$.
However, as is well-known, the image of $\delta$
is generated by the cycle classes of the irreducible components of $V-V^o$.
This shows that $x$ is a linear combination of the cycle classes 
of $e(C)$ and $D$. Since $\ker(d)$ is generated by the cycle classes 
of $e(C)$ and $D$ as $K$-module, so is $E^{11}_2$.
\end{pf}
Let $\langle e(C),f^{-1}(s)\rangle_{s\in Q}\subset H^2_\dR(V)$ denotes the $K$-module
generated by the cycle classes 
of $e(C)$ and irreducible components of $f^{-1}(s)$.
Consider the composition of maps
\begin{equation}\label{VDR-1}
\langle e(C),f^{-1}(s)\rangle_{s\in Q}\lra H^2_\dR(V/K) \lra \prod_{s\in Q}H^2_\dR(f^{-1}(s)).
\end{equation}
This is given by intersection pairing.
Then it is not hard to show that 
\eqref{VDR-1} is injective.
Moreover since the composition
\[
\vg(V,\Omega^2_{V})\lra H^2_\dR(V/K) \lra \prod_{s\in Q}H^2_\dR(f^{-1}(s))
\]
is obviously zero, the second arrow in \eqref{VDR-1}
factors through $E_2^{11}=E_\infty^{11}$.
Summing up this and Lem. \ref{VDR-lem1}, we have a commutative diagram 
\begin{equation}\label{VDR-2}
\xymatrix{
0\ar[r]&\Image\vg(V,\Omega^2_{V})
\ar[r]&H^2_\dR(V/K)\ar[r]& E^{11}_\infty\ar[d]\ar[r]&0&\mbox{(exact)}\\
&&\langle e(C),f^{-1}(s)\rangle_K\ar[r]^{\subset\quad}\ar[u]^\cup\ar[ru]^\cong&
\prod_{s\in Q}  H^2_\dR(f^{-1}(s))
}
\end{equation}
with an exact row.
This shows (1).

\medskip

Next we show (2).
We first prove it in case $V_1\ne X$ and $Q_2\subset S$.
Consider a commutative diagram
\[
\xymatrix{
&&0\ar[d]&0\ar[d]\\
&&H^2_\dR(V_1)_0\ar[d]\ar[r]^{a_1}&H^2_\dR(V_2)_0\ar[d]\\
0\ar[r]&\langle f^{-1}(s)\rangle_{s\in Q_1-Q_2}\ar[r]^{b}&H^2_\dR(V_1)\ar[r]^{a_2}\ar[d]^c
&H^2_\dR(V_2)\ar[r]\ar[d]
&\bigoplus_{s\in Q_1-Q_2}H_{1,\dR}(f^{-1}(s))\\
&&\prod_{s\in Q_1}H^2_\dR(f^{-1}(s))\ar[r]^{\qquad a_3}&H^2_\dR(E_s)
}
\]
where $E_s\subset V_2$ is a fixed smooth fiber and $a_3$ is a projection.
As we have seen in the proof of (1), the composition $cb$ is injective.
Moreover $\Image(c)\cong E^{11}_2$ is generated by
the image of the cycle classes 
of $e(C)$ and irreducible components of $f^{-1}(s)$ with $s\in Q_1$
(Lem. \ref{VDR-lem1}).
Therefore codimension of $\Image(cb)$ in $\Image(c)$ is at most one.
The kernel of $\Image(c)\to H^2_\dR(E_s)$ is of codimension 1 because the cycle
class $[e(C)]$
goes to non-zero via $a_3$.
Hence we have $\Image(c)\cap\ker(a_3)=\Image(cb)$.
Now (2) follows from the snake lemma.
In case $V_1\ne X$ and $V_1\supset V_2$ is arbitrary, we take $V_3\subset V_2\cap U$.
Then a diagram chase of a commutative diagram
\[
\xymatrix{
0\ar[r]& H^2_\dR(V_1)_0\ar[d]\ar[r]&H^2_\dR(V_3)_0\ar[r]\ar@{=}[d]
&\bigoplus_{s\in Q_1-Q_3}H_{1,\dR}(f^{-1}(s))\ar[d]\\
0\ar[r]& H^2_\dR(V_2)_0\ar[r]&H^2_\dR(V_3)_0\ar[r]&
\bigoplus_{s\in Q_2-Q_3}H_{1,\dR}(f^{-1}(s))
}
\]
yields the assertion.
There remains the case $V_1=X$.
However it is easy to see that there is an exact sequence
\[
0\lra H^2_\dR(X)_0\lra H^2_\dR(X-E_s)_0\lra H_{1,\dR}(E_s)
\]
where $E_s=f^{-1}(s)$ is a smooth fiber.
Then the rest of the argument is similar to the above. 

\medskip

Finally we show (3). 
Consider a commutative diagram
\begin{equation}\label{VDR-2-1}
\xymatrix{
0\ar[r]&H^2_\dR(\ol{U})_0\ar[r]&H^2_\dR(\ol{U})\ar[r]&  H^2_\dR(D)\\
0\ar[r]&\NF_\dR(X)^\perp\ar[r]\ar[u]&H^2_\dR(X)\ar[u]^a\ar[r]& \NF_\dR(X)
\ar[u]^b\ar[r]&0\\
}
\end{equation}
with exact rows. 
Since $X-\ol{U}=E$ are additive fibers, $a$ is surjective.
Therefore it is enough to show that $\ker(a)\to \ker(b)$ is bijective. 
$\ker(a)$ is the sub $K$-module generated by the irreducible components
of $E$. This implies $\NF_\dR(X)^\perp\cap \ker(a)=0$ and hence 
$\ker(a)\to \ker(b)$ is injective.
On the other hand, since $\ker(b)$ is generated by
the irreducible components of $E$, $\ker(a)\to \ker(b)$ is surjective.
This completes the proof of (3).
\end{pf}


\subsection{Hodge filtration}\label{hodge-sect}
By taking the embedded resolution of singularities if necessary, 
we can assume that $E_{\mathrm{red}}$ is a NCD.
We then consider the de Rham cohomology groups
\[
H^q_\dR(U)=H^q_\zar(U,\Omega^\bullet_{U})
\cong
H^q_\zar(X,\Omega^\bullet_X(\log D+E))
\]
with the Hodge filtration 
\[
F^p H^q_\dR(U):=\Image[H^q(X,\Omega^{\bullet\geq p}_X(\log D+E))\hra
H^q(X,\Omega^\bullet_X(\log D+E))].
\]
Let $T:=T_m+T_a$.
Define a sheaf $\Omega^1_{X/C}(\log D+E)$ by the exact sequence
\[
0\lra 
f^*\Omega^1_{C}(\log T)\lra
\Omega^1_{X}(\log D+E)\lra
\Omega^1_{X/C}(\log D+E)\lra 0.
\]
This is a locally free sheaf of rank one.
Put
\[
\cH_e:=R^1f_*\Omega^\bullet_{X/C}(\log D+E),\quad
\cH_e^{1,0}:=f_*\Omega^1_{X/C}(\log D+E),\quad
\cH_e^{0,1}:=R^1f_*\O_X.
\]
Then the {\it Gauss-Manin connection}
\[
\nabla:\cH_e\lra \Omega^1_C(\log T)\ot\cH_e
\]
is defined to be the connecting homomorphism arising from an exact sequence  
\begin{equation}\label{hodge-exmp-12}
0\to
f^*\Omega^1_{C}(\log T)\ot \Omega^{\bullet-1}_{X/C}(\log D+E)\to
\Omega^\bullet_{X}(\log D+E)\to
\Omega^\bullet_{X/C}(\log D+E)\to 0
\end{equation}
(see Appendix for a remark on sign.)
Write
\[
H^q_\dR(C,\cH_e):=H^q_\zar(C,\cH_e\to \Omega^1_C(\log T)\ot\cH_e).
\]
\begin{thm}[cf. \cite{SZ} \S 5]\label{Hodge}
Let us put $H^2_\dR(U)_0:=\ker[H^2_\dR(U)\to H^2_\dR(E_s)]$ where $E_s=f^{-1}(s)$
is a smooth fiber contained in $U$.
Then there is the natural isomorphism
\begin{equation}\label{hodge-exmp-8}
H^1_\dR(C,\cH_e)\os{\cong}{\lra}
H^2_\dR(U)_0.
\end{equation}
Moreover under the above isomorphism,
the Hodge filtration corresponds in the following way. 
\begin{align}
&F^1H^2_\dR(U)_0\cong 
H^1_\zar(C,\cH^{1,0}_e
\to \Omega^1_C(\log T)\ot\cH^{1,0}_e)
\rangle\label{hodge-exmp-9}\\
&F^2H^2_\dR(U)_0\cong H^0_\zar(C,\Omega^1_C(\log T)\ot\cH^{1,0}_e)\label{hodge-exmp-10}\\
&\mathrm{Gr}_F^0H^2_\dR(U)_0\cong H^1_\zar(C,\cH^{0,1}_e)\label{hodge-exmp-11}
\end{align}
\end{thm}
\begin{pf}
The exact sequence
\eqref{hodge-exmp-12}
gives rise to a spectral sequence
\[
E_2^{pq}=H^p_\dR(C,R^qf_*\Omega^{\bullet}_{X/C}(\log D))\Longrightarrow
H^{p+q}_\dR(U).
\]
This yields 
\[
0\lra H^1_\dR(C,\cH_e)
\lra H^2_\dR(U)\lra H^0_\dR(C,R^2f_*\Omega^{\bullet}_{X/C}(\log D+E)) 
\lra 0.
\]
Since the last term is one-dimensional, isomorphic to $H^1_\dR(E_s)$, 
we have \eqref{hodge-exmp-8}.

\eqref{hodge-exmp-12} induces an exact sequence 
\[
0\to f^*\Omega^1_C(\log T)\ot
\Omega^{\bullet-1\geq p-1}_{X/C}(\log D+E)
\to 
\Omega^{\bullet\geq p}_X(\log D+E)\lra 
\Omega^{\bullet\geq p}_{X/C}(\log D+E)
\to 0
\]
and this yields
\[
\xymatrix{
H^1_\zar(C,R^1f_*\omega_{X/C}^{\bullet\geq p}
\to \Omega^1_C(\log T)\ot R^1f_*\omega^{\bullet\geq p-1}_{X/C})
\ar[d]\\
H^2_\zar(X,\omega^{\bullet\geq p}_{X/C}
\to f^*\Omega^1_C(\log T)\ot\omega^{\bullet\geq p-1}_{X/C})
\ar[r]^{\qquad\quad\cong}& H^2_\zar(X,\Omega^{\bullet\geq p}_X(\log D+E))\\
}
\]
where $\omega_{X/C}^\bullet:=\Omega^\bullet_{X/C}(\log D+E)$.
Now \eqref{hodge-exmp-9}, 
\eqref{hodge-exmp-10} and \eqref{hodge-exmp-11} easily follow from this.
\end{pf}
A basis of the locally free sheaf $\cH_e$ is given in the following way.
Let $s\in C(\ol{K})$.
We choose a minimal Weierstrass equation
\[
y^2=4x^3-g_2x-g_3,\quad \Delta:=g_2^3-27g_3^2
\]
of $X$ around a (sufficiently small) neighborhood of a fiber $f^{-1}(s)$.
Let $\canh$ and $\can$ be the following elements of $\O_{C,s}\ot\cH_e$
(see \eqref{A1-1} and \eqref{A1-2} in Appendix for the notation):
\[
\canh:=(0)\times (\frac{dx}{y},\frac{dx}{y}),\quad
\can:=(\frac{x^2}{y})\times (\frac{xdx}{2y},\frac{(2g_2x^2+3g_3x)dx}{2y^3}).
\]
Then
\begin{itemize}
\item
If $f^{-1}(s)$ is smooth or multiplicative, then $\{\canh,\can\}$ is a free basis
of $\O_{C,s}\ot\cH_e$.
\item
If $f^{-1}(s)$ is additive, then $\{t\canh,\can\}$ is a basis where $t\in \O_{C,s}$
is a uniformizer.
\end{itemize}
The following theorem is useful.
\begin{thm}[Canonical bundle formula]\label{hodge-exmp}
Let
\begin{center}
\begin{tabular}{c|ccccccccc}
$\epsilon_s$&0&$b$&2&3&4&$b+6$&10&9&8\\
\hline
$f^{-1}(s)$&smooth&$\mathrm{I}_b$&II&III&IV&$\mathrm{I}^*_b$&II*&III*&IV*
\end{tabular}
\end{center}
and put
\begin{equation}\label{hodge-exmp-1}
\epsilon:=\frac{1}{12}\sum_{s\in C(\ol{K})} \epsilon_s\in\Z.
\end{equation}
Then there is an invertible sheaf $\cL$ on $C$ of degree $\epsilon$ such that
\begin{equation}\label{hodge-exmp-1-1}
K_X\cong f^*(K_C\ot\cL),\quad
R^1f_*\O_X\cong\cL^{-1}.
\end{equation}
Mo rover let $a$ be the number of additive fibers 
in the fibration
$f:X_{\ol{K}}\to C_{\ol{K}}$.
Then one has
\begin{equation}\label{hodge-exmp-2}
\deg(\cH^{1,0}_e)=\epsilon-a,\quad \deg(\cH^{0,1}_e)=-\epsilon.
\end{equation} 
\end{thm}

\subsection{Relative cohomology and Extra terms}\label{relative-sect}
For a smooth manifold $M$, we denote by 
$\cA^q(M)$ the space of smooth differential $q$-forms on $M$
with coefficients in $\C$.

\medskip

Suppose $K=\C$.
Let $D_0$ be a union of some multiplicative fibers.
Let $\rho:\wt{D}_0\to D_0$ be the normalization and $\Sigma\subset D_0$
the set of singular points. Let $s:\wt{\Sigma}:=\rho^{-1}(\Sigma)
\hra\wt{D}_0$.
There is the exact sequence
\[
0\lra \O_{D_0}\os{\rho^*}{\lra} \O_{\wt{D}_0}\os{s^*}{\lra} 
\C_{\wt{\Sigma}}/\C_\Sigma\lra 0
\]
where $\C_{\wt{\Sigma}}=\mathrm{Maps}(\wt{\Sigma},\C)=\Hom(\Z\wt{\Sigma},\C)$ etc. and
$\rho^*$ and $s^*$ are the pull-back.
We define $\cA^\bullet(D_0)$ to be the mapping fiber of 
$s^*:\cA^\bullet(\wt{D}_0)\to \C_{\wt{\Sigma}}/\C_\Sigma$:
\[
\cA^0(\wt{D}_0)
\os{s^*\op d}{\lra} \C_{\wt{\Sigma}}/\C_\Sigma\op \cA^1(\wt{D}_0)
\os{0\op d}{\lra} \cA^2(\wt{D}_0)
\]
where the first term is placed in degree 0.
Then
\[
H^q_\dR(D_0)=H^q(\cA^\bullet(D_0))
\]
is the de Rham cohomology of $D_0$, which fits into the exact sequence
\[
\cdots\lra H^0_\dR(\wt{D}_0)\lra \C_{\wt{\Sigma}}/\C_\Sigma
\lra H^1_\dR(D_0)\lra H^1_\dR(\wt{D}_0)\lra\cdots.
\]
There is the natural pairing 
\begin{equation}\label{pairingD}
H_1(D_0,\Z)\otimes H^1_\dR(D_0)\lra \C,
\quad
\gamma\ot z\mapsto\int_\gamma z:=\int_\gamma\eta-c(\partial(\rho^{-1}\gamma))
\end{equation}
where $z=(c,\eta)\in \C_{\wt{\Sigma}}/\C_\Sigma\op\cA^1(\wt{D}_0)$ with $d\eta=0$
and $\partial$ denotes the boundary of homology cycles.

\medskip

Let $V\subset X$ be a Zariski open set containing $D_0$.
We define $\cA^\bullet(V,D_0)$
to be the mapping fiber of 
$j^*:\cA^\bullet(V)\to \cA^\bullet(D_0)$ the pull-back of $j:D_0\hra V$:
\[
\cA^0(V)\os{\cD_0}{\lra} \cA^0(\wt{D}_0)\op \cA^1(V)
\os{\cD_1}{\lra} \C_{\wt{\Sigma}}/\C_\Sigma\op\cA^1(\wt{D}_0)\op \cA^2(V)\os{\cD_2}{\lra}\cdots
\]
where 
\[
\cD_0=j^*\op d,\quad
\cD_1=\begin{pmatrix}
-(s^*\op d)&j^*\\
&d
\end{pmatrix},\quad
\cD_2=\begin{pmatrix}
-(0\op d)&j^*\\
&d
\end{pmatrix},\ldots
\]
Then 
\[
H^q_\dR(V,D_0)=H^q(\cA^\bullet(V,D_0))
\]
is the de Rham cohomology which
fits into the exact sequence
\begin{equation}\label{exp-9}
\cdots\lra H^{q-1}_\dR(D_0)\lra H^q_\dR(V,D_0)\lra H^q_\dR(V)\lra H^q_\dR(D_0)\lra\cdots. 
\end{equation}
In particular, an element of $H^2_\dR(V,D_0)$ is described by $z=(c,\eta,\omega)
\in \C_{\wt{\Sigma}}/\C_\Sigma\op\cA^1(\wt{D}_0)\op \cA^2(V)$ with
$j^*\omega=d\eta$ and $d\omega=0$ which are subject to
relations $(s^*f, df,0)=0$ and $(0,j^*\theta,d\theta)=0$ for $f\in \cA^0(\wt{D}_0)$
and $\theta\in \cA^1(V)$.
The natural pairing
\begin{equation}\label{pairingVD1}
H_2(V,D_0;\C)\ot H_\dR^2(V,D_0)\lra \C,\quad \Gamma\ot z
\longmapsto \int_\Gamma z
\end{equation}
is given by
\begin{equation}\label{pairingVD2}
\int_\Gamma z:=\int_\Gamma \omega-\int_{\partial\Gamma}(c,\eta)
=\int_\Gamma \omega-\int_{\partial\Gamma}\eta+c(\rho^{-1}(\partial\Gamma)).
\end{equation}

\medskip

There are canonical maps
\begin{equation}\label{pairingVD7}
\vg(V,\Omega^2_V)\lra H^2_\dR(V)=H^2(\cA^\bullet(V)),\quad \omega\longmapsto 
\ol{\omega}.
\end{equation}
\begin{equation}\label{pairingVD5}
\vg(V,\Omega^2_V)\lra H^2_\dR(V,D_0),\quad \omega\longmapsto 
(0,0,\omega).
\end{equation}
Define $G\vg(V,\Omega^2_V)\subset \vg(V,\Omega^2_V)$ to be the inverse image of $F^1H^2_\dR(V)$
via the natural map $\vg(V,\Omega^2_V)\to H^2_\dR(V)$ 
where $F^\bullet$ denotes the Hodge filtration.
We define a map $\cand_{D_0}$ by  
a commutative diagram
\begin{equation}\label{pairingVD3}
\xymatrix{
0\ar[r]&G\vg(V,\Omega^2_V)\ar[d]^{\cand_{D_0}}\ar[r]&\vg(V,\Omega^2)\ar[d]^{\eqref{pairingVD5}}
\ar[r]&\Gr^0_FH^2_\dR(V)\ar@{=}[d]\\
0\ar[r]&H^1_\dR(D_0)\ar[r]^{i^0\quad}&
\Gr^0_FH^2_\dR(V,D_0) \ar[r]&\Gr^0_F H^2_\dR(V)\ar[r]&0\\
}
\end{equation}
with exact rows. Here $i:H^1_\dR(D_0)\to H^2_\dR(V,D_0)$ is the map appearing in
\eqref{exp-9} and $i^0$ denotes the induced map on the graded piece.
We call $\cand_{D_0}(\omega)$ the {\it extra term} of $\omega$
at $D_0$.
\begin{prop}\label{pairingVD4}
Let $\omega\in G\vg(V,\Omega^2_V)$.
Then $\omega_{V,D_0}:=(0,0,\omega)-i\cand_{D_0}(\omega)\in F^1H^2_\dR(V,D_0)$ is the unique 
element corresponding to $\ol{\omega}\in H^2_\dR(V)$ 
via the natural map $F^1H^2_\dR(V,D_0)\to F^1H^2_\dR(V)$.
Moreover
\begin{equation}\label{pairingVD6}
\int_\Gamma\omega_{V,D_0}=\int_\Gamma \omega-\int_{\partial\Gamma}\cand_{D_0}(\omega).
\end{equation}
\end{prop}
\begin{pf}
It follows from the construction that 
$\omega-i\cand_{D_0}(\omega)$ belongs to $F^1H^2_\dR(V,D_0)$.
The uniqueness follows from the injectivity of the map
$F^1H^2_\dR(V,D_0)\to F^1H^2_\dR(V)$.
\eqref{pairingVD6} follows from \eqref{pairingVD2}.
\end{pf}
The map ``$\cand_{D_0}$" can be defined in an algebraic way. 
Let us denote by $(\check{C}^\bullet({\mathscr F}),\delta)$ 
the Cech complex of a sheaf $\mathscr F$.
Then $H^1_\dR(D_0)$ is isomorphic to the cohomology of the complex
\[
\check{C}^0(\O_{\wt{D}_0})
\os{\cD_0}{\lra}
\check{C}^1(\O_{\wt{D}_0})
\times
\check{C}^0(\C_{\wt{\Sigma}}/\C_\Sigma\op \Omega^1_{\wt{D}_0})
\os{\cD_1}{\lra}
\check{C}^2(\O_{\wt{D}_0})
\times
\check{C}^1(\C_{\wt{\Sigma}}/\C_\Sigma\op \Omega^1_{\wt{D}_0})
\]
at the middle term where
\[
\cD_0=\delta\times(s^*\op d),\quad
\cD_1:=\begin{pmatrix}
\delta&-(s^*\op d)\\
&\delta
\end{pmatrix}.
\]
Moreover $H^2_\dR(V)$ and $H^2_\dR(V,D_0)$ are isomorphic to the cohomology
of the following complexes
\[
\check{C}^1(\O_V)\times
\check{C}^0(\Omega^1_V)
\os{\cD_2}{\lra}
\check{C}^2(\O_V)\times
\check{C}^1(\Omega^1_V)
\times
\check{C}^0(\Omega^2_V)
\os{\cD_3}{\lra}
\check{C}^3(\O_V)\times
\check{C}^2(\Omega^1_V)
\times
\check{C}^1(\Omega^2_V)
\]
\begin{multline*}
\check{C}^1(\O_V)\times
\check{C}^0(\O_{\wt{D}_0}\op\Omega^1_V)
\os{\cD_4}{\lra}
\check{C}^2(\O_V)\times
\check{C}^1(\O_{\wt{D}_0}\op\Omega^1_V)
\times
\check{C}^0(\C_{\wt{\Sigma}}/\C_\Sigma\op\Omega^1_{\wt{D}_0}\op\Omega^2_V)\\
\os{\cD_5}{\lra}
\check{C}^3(\O_V)\times
\check{C}^2(\O_{\wt{D}_0}\op\Omega^1_V)
\times
\check{C}^1(\C_{\wt{\Sigma}}/\C_\Sigma\op\Omega^1_{\wt{D}_0}\op\Omega^2_V)
\end{multline*}
at the middle terms respectively,
\begin{equation}\label{d4}
\cD_2=
\begin{pmatrix}
\delta&-d&\\
&\delta&d
\end{pmatrix}
,\quad
\cD_3=
\begin{pmatrix}
\delta&d&\\
&\delta&-d\\
&&\delta
\end{pmatrix}
,\quad\cD_4=
\begin{pmatrix}
\delta&-(j^*\op d)&\\
&\delta&T
\end{pmatrix},
\end{equation}
\begin{equation}\label{d5}
\cD_5=
\begin{pmatrix}
\delta&j^*\op d&\\
&\delta&-T\\
&&\delta
\end{pmatrix},\quad
T=
\begin{pmatrix}
-s^*&-d&\\
&j^*&d
\end{pmatrix}.
\end{equation}
For a $\omega\in G\vg(V,\Omega^2_V)$, we simply write 
$\omega=(0)\times(0)\times(\omega)\in 
\check{C}^2(\O_V)\times \check{C}^1(\Omega^1_V)\times \check{C}^0(\Omega^2_V)$.
There is a Cech cocycle
\[
\xi=(0)\times (\eta_{ij})\times (\omega_i)
\in 
\check{C}^2(\O_X)\times \check{C}^1(\Omega^1_X)\times \check{C}^0(\Omega^2_X)
\]
such that $\xi\equiv
\omega$ in $H^2_\dR(V)$
and $\xi$ belongs to $F^1H^2_\dR(X)$ and the kernel of $H^2_\dR(X)\to H^2_\dR(D_0)$
(Prop \ref{VDR} (1)).
Then there is a unique Cech cocycle
\[
\xi_{X,D_0}=
(0)\times (0,\eta_{ij})\times (0,\wt{\eta}_i,\omega_i)
\in \check{C}^2(\O_X)\times \check{C}^1(\O_{\wt{D}_0}\op\Omega^1_X)
\times \check{C}^0(\C_{\wt{\Sigma}}/\C_\Sigma\op\Omega^1_{\wt{D}_0}\op\Omega^2_X)
\]
such that $\eta_{ij}|_{\wt{D}_0}=\wt{\eta}_j-\wt{\eta}_i$.
This belongs to $F^1H^2_\dR(X,D_0)$ by definition.
Since $\xi\equiv\omega$ in $H^2_\dR(V)$,
there is a Cech cycle $z=(f_{ij})\times (\nu_i)
\in \check{C}^1(\O_V)\times \check{C}^0(\Omega^1_V)$ such that
\begin{align*}
\omega-\xi
&=(0)\times (-\eta_{ij})\times (\omega-\omega_i)\\
&=(f_{jk}-f_{ik}+f_{ij})\times (-df_{ij}+\nu_j-\nu_i)\times (d\nu_i)
\end{align*}
in $\check{C}^2(\O_V)\times \check{C}^1(\Omega^1_V)\times \check{C}^0(\Omega^2_V)$.
Hence
\begin{align*}
(0)\times(0,0)\times(0,0,\omega)-\xi_{X,D_0}
&=(0)\times (0,-df_{ij}+\nu_j-\nu_i)\times (0,-\wt{\eta}_{i},d\nu_i)\\
&\equiv(0)\times (f_{ij}|_{\wt{D}_0},0)\times (0,-\wt{\eta}_{i}-\nu_i|_{\wt{D}_0},0)\\
&\in
\check{C}^2(\O_V)\times \check{C}^1(\O_{\wt{D}_0}\op
\Omega^1_V)\times \check{C}^0(\C_{\wt{\Sigma}}/\C_\Sigma\op\Omega^1_{\wt{D}_0}\op\Omega^2_V)
\end{align*}
modulo the image of $\check{C}^1(\O_V)\times \check{C}^0(\O_{\wt{D}_0}\op\Omega^1_V)$.
This shows
$
\cand_{D_0}(\omega)=(f_{ij}|_{\wt{D}_0})\times (0,-\wt{\eta}_{i}-\nu_i|_{\wt{D}_0})
$
in $H^1_\dR(D_0)$.
There is $(h_i)\in \check{C}^0(\O_{\wt{D}_0})$ such that
$f_{ij}|_{\wt{D}_0}=h_j-h_i$.
Then 
\[
\cand_{D_0}(\omega)=(f_{ij}|_{\wt{D}_0})\times (0,-\wt{\eta}_{i}-\nu_i|_{\wt{D}_0})
\equiv(0)\times (s^*h_i,0)
\in H^1_\dR(D_0).
\]

\def\goodua{\Lambda^1({U})}
\def\goodub{\Lambda^2(U)}
\def\goodxa{\Lambda^1(\ol{U})}
\def\goodxb{\Lambda^2(\ol{U})}
\subsection{Good algebraic 2-forms}\label{good-sect}
We introduce two subspaces
\[
\goodub\subset\goodua\subset\vg(S^o,\Omega^1_{S^o}\ot\cH_e)=\vg(U^o,\Omega^2_{U^o}),
\]
which we call the spaces of {\it good algebraic 2-forms}.
Define
\[
\goodub:=\Image[\vg(C,\Omega^1_{C}(\log T)\ot\cH_e^{1,0})
=\vg(X,\Omega^2_X(\log D+E))\hra \vg(U^o,\Omega^2_{U^o})].
\]
We define $\goodua$ in the following way.
Let us consider a diagram
\begin{equation}\label{hodge-pf1}
\xymatrix{
&0\ar[d]\\
&\Omega^1_{S^o}\ot\cH_e^{1,0}|_{S^o}\ar[d]\\
\cH_e^{1,0}|_{S^o}\ar[r]^{\nabla\quad}\ar[d]_=&\Omega^1_{S^o}\ot\cH_e|_{S^o}\ar[d]\\
\cH_e^{1,0}|_{S^o}\ar[r]^{\ol{\nabla}\quad}&\Omega^1_{S^o}\ot\cH_e^{0,1}|_{S^o}\ar[d]\\
&0}
\end{equation}
It follows from the definition of $S^o$ and Cor. \ref{PF-rem}
that the bottom arrow $\ol{\nabla}$ is isomorphism.
This yields an isomorphism
\begin{equation}\label{hodge-pf2}
\Omega^1_{S^o}\ot\cH_e^{1,0}|_{S^o}\os{\cong}{\lra}
\Omega^1_{S^o}\ot\cH|_{S^o}/\Image(\cH^{1,0}_e|_{S^o})
\end{equation}
and 
\begin{align*}
H^1_\zar(C,\cH_e^{1,0}\to\Omega^1_C(\log T)\ot\cH_e)
&\lra
H^1_\zar(S^o,\cH_e^{1,0}\to\Omega^1_{S^o}\ot\cH_e)\\
&\os{\cong}{\lra}
\vg(S^o,\Omega^1_{S^o}\ot\cH_e/\Image(\cH^{1,0}_e|_{S^o}))\\
&\us{\cong}{\os{\eqref{hodge-pf2}}{\longleftarrow}}
\vg(S^o,\Omega^1_{S^o}\ot\cH_e^{1,0})=\vg(U^o,\Omega^2_{U^o}).
\end{align*}
Define $\goodua$ to be the image of the composition of the above maps:
\[
\goodua:=\Image[H^1_\zar(C,\cH_e^{1,0}\to\Omega^1_C(\log T)\ot\cH_e)
\to \vg(U^o,\Omega^2_{U^o})].
\]
\begin{prop}\label{good-prop1}
\[
H^1_\zar(C,\cH_e^{1,0}\to\Omega^1_C(\log T)\ot\cH_e)
\os{\cong}{\lra}\goodua,\quad
\vg(X,\Omega^2_X(\log D+E))
\os{\cong}{\lra}\goodub.
\]
Hence
\[
\goodua\os{\cong}{\lra} F^1H^2_\dR(U)_0,\quad
\goodub\os{\cong}{\lra} F^2H^2_\dR(U)_0
\]
by Thm \ref{Hodge}.
\end{prop}
\begin{pf}
There is nothing to show other than the injectivity of
$H^1_\zar(C,\cH_e^{1,0}\to\Omega^1_C(\log T)\ot\cH_e)
\to\goodua
$. However this follows from the fact that $F^1H^2_\dR(U)\to F^1H^2(U^o)$ is 
injective.
\end{pf}

\begin{lem}\label{hodge-pf-lem1}
Along $D$, good algebraic 2-forms have at most log poles.
Namely
\begin{equation}\label{hodge-pf3}
\goodua\subset 
\vg(\ol{U^o},\Omega^2_{\ol{U^o}}(\log D)).
\end{equation}
\end{lem}
\begin{pf}
We may replace $K$ with $\ol{K}$.
Let us consider a diagram
\begin{equation}\label{hodge-pf1-d}
\xymatrix{
&0\ar[d]\\
&\Omega^1_{\ol{S^o}}(\log T_m)\ot\cH_e^{1,0}|_{\ol{S^o}}\ar[d]\\
\cH_e^{1,0}|_{\ol{S^o}}\ar[r]^{\nabla\qquad\quad}\ar[d]_=
&\Omega^1_{\ol{S^o}}(\log T_m)\ot\cH_e|_{\ol{S^o}}\ar[d]\\
\cH_e^{1,0}|_{\ol{S^o}}\ar[r]^{\ol{\nabla}\qquad\quad}&
\Omega^1_{\ol{S^o}}(\log T_m)\ot\cH_e^{0,1}|_{\ol{S^o}}\ar[d]\\
&0}
\end{equation}
Around a multiplicative fiber $D_0=f^{-1}(s_0)$, $X$ can be written by a Weierstrass form
$y^2=4x^3-g_2x-g_3$ with $\ord_{s_0}(g_2^3-27g_3^2)>0$ and $\ord_{s_0}(g_2)=
\ord_{s_0}(g_3)=0$ where $\ord_{s_0}$ denotes the valuation order
on $\O_{C,s_0}$ (cf. Tate's algorithm).
Thus Cor. \ref{PF-rem} implies that the bottom arrow $\ol{\nabla}$ in \eqref{hodge-pf1-d}
is an isomorphism.
Then we have  
\begin{align*}
H^1_\zar(C,\cH_e^{1,0}\to\Omega^1_C(\log T)\ot\cH_e)
&\lra
H^1_\zar(\ol{S^o},\cH_e^{1,0}\to\Omega^1_{\ol{S^o}}(\log T_m)\ot\cH_e)\\
&\os{\cong}{\lra}
\vg(\ol{S^o},\Omega^1_{\ol{S^o}}(\log T_m)\ot\cH_e/\Image(\cH^{1,0}_e|_{\ol{S^o}}))\\
&\os{\cong}{\longleftarrow}
\vg(\ol{S^o},\Omega^1_{\ol{S^o}}(\log T_m)\ot\cH_e^{1,0})\\
&=
\vg(\ol{U^o},\Omega^2_{\ol{U^o}}(\log D))\\
&\hra\vg(U^o,\Omega^2_{U^o})
\end{align*}
and this shows \eqref{hodge-pf3}.
\end{pf}
By Lem.\ref{hodge-pf-lem1}, 
one can have the residue map
\begin{equation}\label{hodge-pf4}
\Res_D:\goodua\lra H_{1,\dR}(D)
\end{equation} 
along $D$.
We define
\begin{equation}\label{hodge-pf5}
\goodxa:=\goodua\cap\ker(\Res_D)\supset
\goodxb:=\goodub\cap\ker(\Res_D).
\end{equation} 
\begin{lem}\label{hodge-pf-lem2}
$\goodxa\subset \vg(\ol{U^o},\Omega^2_{\ol{U^o}})$.
\end{lem}
\begin{pf}
There is the weight filtration $W_\bullet\Omega^2_{\ol{U^o}}(\log D)$ such that
$\Gr^W_0\Omega^2_{\ol{U^o}}(\log D)=\Omega^2_{\ol{U^o}}$,
$\Gr^W_1\Omega^2_{\ol{U^o}}(\log D)=\Omega^1_{\wt{D}}$ and
$\Gr^W_2\Omega^2_{\ol{U^o}}(\log D)=\Omega^2_{\wt{\Sigma}}$.
By definition one has $\goodxa\subset \vg(\ol{U^o},W_1\Omega^2_{\ol{U^o}}(\log D))$.
However since $\wt{D}$ is a union of $\P^1$, one has 
$\vg(\ol{U^o},\Gr_1^W\Omega^2_{\ol{U^o}}(\log D))=\vg(\wt{D},\Omega^1_{\wt{D}})=0$
and hence
$\vg(\ol{U^o},W_1\Omega^2_{\ol{U^o}}(\log D))=
\vg(\ol{U^o},W_0\Omega^2_{\ol{U^o}}(\log D))
=\vg(\ol{U^o},\Omega^2_{\ol{U^o}})$.
\end{pf}
By Prop. \ref{VDR} (1), the image of $\goodxa$ via the natural map
$\vg(\ol{U^o},\Omega^2_{\ol{U^o}})
\to H^2_\dR(\ol{U^o})$ is contained in $H^2_\dR(\ol{U^o})_0$.
\begin{prop}\label{good-prop2}
\[\goodxa \cong F^1H^2_\dR(\ol{U})_0,\quad
\goodxb\cong F^2H^2_\dR(\ol{U})_0.\]
\end{prop}
\begin{pf}
Prop. \ref{VDR} (2) and 
the definition of $\goodxa$ give rise to
a commutative diagram
\[
\xymatrix{
0\ar[r]&\goodxa\ar[d]\ar[r]&\goodua\ar[d] \ar[r] &H_{1,\dR}(D)\ar@{=}[d]\\
0\ar[r]&F^1H^2_\dR (\ol{U})_0 \ar[r]&F^1H^2_\dR (U)_0\ar[r]&H_{1,\dR}(D)
}
\]
with exact rows. Now the assertion follows from 
Prop. \ref{good-prop1}.
\end{pf}
The following proposition is the motivation by which we introduced the good algebraic 2-forms.
\begin{prop}\label{hodge-mot}
Let $\omega\in \goodxa$ be a good 2-form.
Suppose that $D_0=f^{-1}(s_0)$ is an irreducible multiplicative fiber
(i.e. type $I_1$).
Then the extra term $\cand_{D_0}(\omega)$ is zero.
\end{prop}
\begin{rem}\label{hodge-mot-rem}
Prop. \ref{hodge-mot} seems true for a fiber of type $I_n$ for arbitrary $n\geq 1$.
\end{rem}
\begin{pf}
We may replace $K$ with $\ol{K}$.
Put $E^*=E+f^{-1}(F)$.
We use the description of $H^\bullet_\dR(X)$ etc. by the Cech complexes.
Let 
\[
(\alpha_{ij})\times (\beta_i)\in \check{C}^1(\cH^{1,0}_e)\times
\check{C}^0(\Omega^1_C(\log T)\ot\cH_e)
\]
be a corresponding Cech cocycle to $\omega$, and this defines
\[
z:=(0)\times (\eta_{ij})\times (\pi_i)\in \check{C}^2(\O_X)\times
\check{C}^1(\Omega^1_X(\log D+E^*))\times
\check{C}^0(\Omega^2_X(\log D+E^*))
\]
in a natural way.
The proof of Lem. \ref{hodge-pf-lem1} shows that there is 
$y_0\in \check{C}^0(\cH_e|_{\ol{S^o}})$ such that 
\[
(0)\times (\omega)=(\alpha_{ij})\times (\beta_i)-\cD_0(y_0)
\]
where $\cD_0:\check{C}^0(\cH_e)\to
\check{C}^1(\cH^{1,0}_e)\times
\check{C}^0(\Omega^1_C(\log T)\ot\cH_e)$.
This means that there is $y=(0)\times(\nu_i)
\in \check{C}^1(\O_{\ol{U^o}})\times\check{C}^0(\Omega^1_{\ol{U^o}}(\log D))$
such that
\begin{equation}\label{hodge-pf8}
z|_{\ol{U^o}}-\cD(y)=
(0)\times (\eta_{ij}|_{\ol{U^o}}-(\nu_j-\nu_i))\times (\pi_i|_{\ol{U^o}}-d\nu_i)
=(0)\times(0)\times (\omega)
\end{equation}
where
\[
\cD:\check{C}^1(\O_X)\times\check{C}^0(\Omega^1_X(\log D+E^*))
\to
\check{C}^2(\O_X)\times
\check{C}^1(\Omega^1_X(\log D+E^*))\times
\check{C}^0(\Omega^2_X(\log D+E^*)).
\]
On the other hand, there is a Cech cocycle
$
w=(0)\times(*)\times(*)\in \check{C}^2(\O_X)\times
\check{C}^1(\Omega^1_X(\log E^*))\times
\check{C}^0(\Omega^2_X(\log E^*))
$
such that $[w]\in F^1H^2_\dR(\ol{U^o})$ belongs to the kernel of
$H^2_\dR(\ol{U^o})\to H^2_\dR(D)$ and $[w]|_U=[z]$ in $H^2_\dR(U)$.
Since $[w]|_U=[z]$ in 
$F^1H^2_\dR(U)\cong H^2(X,\Omega^{\bullet\geq 1}_X(\log D+E^*))$, this means that 
there is $\wt{y}=(0)\times (\wt{\nu}_i)\in \check{C}^1(\O_X)\times
\check{C}^0(\Omega^1_X(\log D+E^*))$
such that 
\begin{equation}\label{hodge-pf9}
w=z-\cD(\wt{y})=
(0)\times (\eta_{ij}-(\wt{\nu}_j-\wt{\nu}_i))\times (\pi_i-d\wt{\nu}_i)
\end{equation}
and this belongs to $\check{C}^2(\O_X)\times
\check{C}^1(\Omega^1_X(\log E^*))\times
\check{C}^0(\Omega^2_X(\log E^*))$.
\begin{lem}\label{extra2}
Fix an arbitrary multiplicative fiber $D_0=f^{-1}(s_0)$, and choose a
(sufficiently small) neighborhood $V$ of $D_0$.
Then there is a constant $c$ such that
\[
\theta_i:=\wt{\nu}_i|_V-\nu_i|_V-c\frac{dt}{t-s_0}
\]
has no log pole along $D_0$.
\end{lem}
\begin{pf}
There is the exact sequence
\[
0\lra\Omega^1_V\lra
\Omega^1_V(\log D_0)\os{\Res}{\lra} \O_{\wt{D}_0}\lra 0.
\]
Since neither $z|_{\ol{U^o}}-\cD(y)$ or $z-\cD(\wt{y})$ has log pole along $D_0$, 
\[
\Res(\eta_{ij})=\Res(\nu_j)-\Res(\nu_i)
=\Res(\wt{\nu}_j)-\Res(\wt{\nu}_i)\in \check{C}^1(\O_{\wt{D}_1}).
\]
Since $D_0$ is irreducible, 
$\ker[\check{C}^0(\O_{\wt{D}_0})\to \check{C}^1(\O_{\wt{D}_0})]$ is one-dimensional,
and hence 
$\Res(\nu_i)-\Res(\wt{\nu}_i)$ is a constant $c$. This implies that
\[
\theta_i:=\nu_i|_V-\wt{\nu}_i|_V-c\frac{dt}{t-s_0}
\]
has no log pole.
\end{pf}
We turn to the proof of Prop. \ref{hodge-mot}.
By Lem. \ref{extra2} and \eqref{hodge-pf8} and \eqref{hodge-pf9},
one has
\begin{align*}
z-\cD(\wt{y})|_V&=
(0)\times ((\eta_{ij}-(\nu_j-\nu_i))|_V-(\theta_j-\theta_i))
\times ((\pi_i-d\nu_i)|_V-d\theta_i)\\
&=(0)\times (-(\theta_j-\theta_i))
\times (\omega-d\theta_i).
\end{align*}
Let $z_{X,D_0}\in F^1H^2_\dR(\ol{U^o},D_0)$ be the corresponding Cech
cocycle to $z-\cD(\wt{y})$ via the injection 
$F^1H^2_\dR(\ol{U^o},D_0)\hra F^1 H^2_\dR(\ol{U^o})$. Then
\begin{align*}
z_{X,D_0}|_V&=(0)\times(0,-(\theta_j-\theta_i))
\times
(0,-(\theta_j-\theta_i)|_{\wt{D}},\omega|_V-d\theta_i)\\
&\equiv 
(0)\times(0,0)\times(0,0,\omega|_V)\mbox{ in }H^2_\dR(V,D_0).
\end{align*}
This belongs to $F^1H^2_\dR(V,D_0)$ since so does $z_{X,D_0}$.
Hence $(0)\times(0,0)\times(0,0,\omega|_V)\in F^1H^2_\dR(V,D_0)$ and this means
$\cand_{D_0}(\omega)=0$. 
\end{pf}


\section{Explicit computations of regulator on $K_1$ of elliptic surfaces}\label{exp-sect}
We keep the notations in \S \ref{deRh-sect-1}.
The base field $K$ is $\C$ and
we assume $D\ne\emptyset$ throughout this section.

\subsection{1-Extension of MHS's arising from a multiplicative
fiber}\label{exp-sect1}
For each $\gamma\in H_1(D,\Q)$, there is a corresponding element $\xi_\gamma\in H^3_\cM(X,\Q(2))$
which is unique up to 
the decomposable part. It is given in the following way.
Let $D_k=\sum_{i=1}^n D_k^{(i)}=f^{-1}(P_k)$ be a multiplicative
fiber over a point $P_k\in C$ and $Q_i$ the intersection points.
There are rational functions $f_i$ on $D_k^{(i)}$ such that 
$\Div_{D_k^{(i)}}(f_i)=Q_{i+1}-Q_i$.
Then we put
\[
\xi_{D_k}:=\sum_{i=1}^n[f_i,D_k^{(i)}]\in H^3_\cM(X,\Q(2)).
\]
If $\gamma=(m_1,\cdots,m_s)\in H_1(D,\Z)=\bigoplus_k H_1(D_k,\Z)\cong\Z^{\op k}$
then we put $\xi_\gamma:=m_1\xi_{D_1}+\cdots +m_s\xi_{D_s}$.
This depends on the choice of $f_i$'s, though the ambiguity 
is killed by the decomposable part.

\begin{center}
\unitlength 0.1in
\begin{picture}( 24.0000, 19.2000)( 41.4000,-21.0000)
%
{\color[named]{Black}{%
\special{pn 8}%
\special{pa 5120 370}%
\special{pa 4330 1190}%
\special{fp}%
\special{pa 4380 920}%
\special{pa 4820 1950}%
\special{fp}%
}}%
%
{\color[named]{Black}{%
\special{pn 8}%
\special{pa 5670 370}%
\special{pa 6460 1190}%
\special{fp}%
\special{pa 6410 920}%
\special{pa 5970 1950}%
\special{fp}%
}}%
%
{\color[named]{Black}{%
\special{pn 8}%
\special{pa 4790 480}%
\special{pa 5950 480}%
\special{fp}%
}}%
%
{\color[named]{Black}{%
\special{pn 4}%
\special{sh 1}%
\special{ar 4990 2010 10 10 0  6.28318530717959E+0000}%
\special{sh 1}%
\special{ar 5276 2080 10 10 0  6.28318530717959E+0000}%
\special{sh 1}%
\special{ar 5566 2100 10 10 0  6.28318530717959E+0000}%
\special{sh 1}%
\special{ar 5846 2020 10 10 0  6.28318530717959E+0000}%
}}%
%
{\color[named]{Black}{%
\special{pn 4}%
\special{sh 1}%
\special{ar 4442 1070 16 16 0  6.28318530717959E+0000}%
\special{sh 1}%
\special{ar 4442 1070 16 16 0  6.28318530717959E+0000}%
\special{sh 1}%
\special{ar 4442 1070 16 16 0  6.28318530717959E+0000}%
}}%
%
{\color[named]{Black}{%
\special{pn 4}%
\special{sh 1}%
\special{ar 5014 482 8 8 0  6.28318530717959E+0000}%
\special{sh 1}%
\special{ar 5014 482 8 8 0  6.28318530717959E+0000}%
}}%
%
{\color[named]{Black}{%
\special{pn 4}%
\special{sh 1}%
\special{ar 5016 482 16 16 0  6.28318530717959E+0000}%
\special{sh 1}%
\special{ar 5016 482 16 16 0  6.28318530717959E+0000}%
}}%
%
{\color[named]{Black}{%
\special{pn 4}%
\special{sh 1}%
\special{ar 5774 476 16 16 0  6.28318530717959E+0000}%
\special{sh 1}%
\special{ar 5774 476 16 16 0  6.28318530717959E+0000}%
}}%
%
{\color[named]{Black}{%
\special{pn 4}%
\special{sh 1}%
\special{ar 6344 1070 16 16 0  6.28318530717959E+0000}%
\special{sh 1}%
\special{ar 6344 1070 16 16 0  6.28318530717959E+0000}%
}}%
\put(41.6000,-10.6000){\makebox(0,0)[lb]{$Q_n$}}%
\put(44.5000,-7.6000){\makebox(0,0)[lb]{$D^{(n)}_k$}}%
\put(48.7000,-3.9000){\makebox(0,0)[lb]{$Q_1$}}%
\put(53.0000,-4.5000){\makebox(0,0)[lb]{$D_k^{(1)}$}}%
\put(58.0000,-3.9000){\makebox(0,0)[lb]{$Q_2$}}%
\put(60.3000,-7.2000){\makebox(0,0)[lb]{$D^{(2)}_k$}}%
\put(62.3000,-16.0000){\makebox(0,0)[lb]{}}%
\put(64.7000,-10.7000){\makebox(0,0)[lb]{$Q_3$}}%
\end{picture}%
\end{center}

\bigskip

Let us recall the regulator $\reg(\xi_\gamma)$ 
from \S \ref{bregin-sect}.
Let
\[
\NF^B(X):=\Image(H_2(D,\Q)\op H_2(E,\Q)\op H_2(e(C),\Q)\lra H_2(X,\Q)).
\]
Then there are the natural isomorphisms
\begin{equation}\label{exp-lem0-2}
H_2(X,\Q)/\NF^B(X)\cong (\NF(X)^\perp)^*\cong \NF(X)^\perp\ot\Q(2).
\end{equation}
The exact sequence
\begin{equation}\label{exp-1}
\xymatrix{
0\ar[r]&H_2(X,\Q)/H_2(D,\Q)\ar[r]
&H_2(X,D,\Q)\ar[r]^\partial&H_1(D,\Q)\ar[r]&0}
\end{equation}
of mixed Hodge structures gives rise to a map
\begin{align*}
\rho:H_1(D,\Q)&\lra \Ext^1(\Q,H_2(X,\Q)/H_2(D,\Q))\\
&\lra\Ext^1(\Q,H_2(X,\Q)/\NF^B(X))\\
&\os{\cong}{\longleftarrow}\Ext^1(\Q,\NF(X)^\perp\ot\Q(2))
\end{align*}
where $\NF(X)^\perp\subset H^2(X,\Q)$ is a Hodge structure of weight 2.
 Then we have from Thm.\ref{bregin-thm}
\begin{equation}\label{exp-1-p}
\reg(\xi_\gamma)=\pm\rho(\gamma).
\end{equation}
In this section, we shall use a slight modification of \eqref{exp-1}.
\begin{lem}\label{exp-lem0}
$\partial:H_2(\ol{U},D;\Q)\to H_1(D,\Q)$ is surjective.
\end{lem}
\begin{pf}
The assertion is equivalent to saying that
$H^1(\ol{U},D;\Q)\to H^1(\ol{U},\Q)$ is surjective.
Since the functional $j$-invariant of $U/S$ is not constant (by the assumption),
one has  
\[
\vg(S,R^1f_*\Q)=H^1(f^{-1}(t),\Q)^{\pi_1(S,t)}=0
\]
and hence $H^1(U,\Q)=H^1(S,\Q)$. This and a commutative diagram
\[
\xymatrix{
0=H^1_D(\ol{U},\Q)\ar[r]& H^1(\ol{U},\Q)\ar[r]& H^1(U,\Q)\ar[r]& H^2_D(\ol{U},\Q)
\cong H_2(D,\Q)\\
0=H^1_{T_m}(\ol{S},\Q)\ar[r]\ar[u]& H^1(\ol{S},\Q)\ar[r]\ar[u]& H^1(S,\Q)
\ar[r]\ar[u]^\cong& H^2_{T_m}(\ol{S},\Q)\cong H_0(T_m,\Q)\ar[u]^\cup
}\]
yield 
\begin{equation}\label{exp-lem0-1}
H^1(\ol{S},\Q)\os{\cong}{\lra}H^1(\ol{U},\Q).
\end{equation}
Therefore, to show the surjectivity of $H^1(\ol{U},D;\Q)\to H^1(\ol{U},\Q)$ 
it is enough to show that
$H^1(\ol{S},T_m;\Q)\to H^1(\ol{S},\Q)$ is surjective. 
However it is clear because $H^1(T_m,\Q)=0$.
\end{pf}

\medskip

\begin{align}
\Ext^1(\Q,\NF(X)^\perp\ot\Q(2))&=\Coker[\NF(X)^\perp\ot\Q(2)\lra
\NF(X)^\perp\ot\C/F^2]\\
&\cong\Coker[H_2(X,\Q)/\NF^B(X)\to
\Hom(F^1H^2_\dR(\ol{U})_0,\C)]\label{exp-ext2}\\
&\cong\Coker[H_2(\ol{U^o},\Q)\lra
\Hom(F^1H^2_\dR(\ol{U})_0,\C)]\label{exp-ext1}\\
&\cong\Coker[H_2(\ol{U^o},\Q)\os{\Phi}{\lra}
\Hom(\goodxa,\C)]\label{exp-6}
\end{align}
where \eqref{exp-ext2} follows from Prop.\ref{VDR} (3) and \eqref{exp-lem0-2},
and 
\eqref{exp-ext1} follows from the surjectivity of 
$H_2(\ol{U^o},\Q)\to H_2(X,\Q)/\NF^B(X)$ and
\eqref{exp-6} follows from Prop.\ref{good-prop2}.
The map $\Phi$ in \eqref{exp-6} is given by
\[
\Delta\longmapsto
\left[\omega\longmapsto\int_\Delta \omega\right],\quad \omega\in \goodxa
\footnote{Note that, since $\goodxa\subset
\vg(\ol{U^o},\Omega^2_{\ol{U^o}})$, one can {\it a priori} define ``$\int_\Delta\omega$" 
only for $\Delta\in H_2(\ol{U^o},\Q)$.}.
\]
Next consider a commutative diagram
\begin{equation}\label{main-diagram}
\xymatrix{
0\ar[r]&H_2(X,\Q)/H_2(D,\Q)\ar[d]_{\mbox{surj.}}\ar[r]\ar@{}[rd]|{\square}
&H_2(X,D;\Q)\ar[d]_{\mbox{surj.}}\ar[r]&H_1(D,\Q)\ar@{=}[d]\ar[r]&0\\
0\ar[r]&H_2(X,\Q)/\NF^B(X)\ar[r]&\ol{H_2(X,D;\Q)}\ar[r]&H_1(D,\Q)\ar[r]&0\\
0\ar[r]&H_2(\ol{U^o},\Q)/H_2(D)\ar[u]^{\mbox{surj.}}\ar[r]
&H_2(\ol{U^o},D;\Q)\ar[u]^{\mbox{surj.}}\ar[r]^\partial&
H_1(D,\Q)\ar[r]\ar@{=}[u]&0
}
\end{equation}
where the surjectivity of the right arrows follows from Lem. \ref{exp-lem0}.
The middle row gives the regulator class \eqref{exp-1-p}.
Let us describe it explicitly under the identification \eqref{exp-6}.
There is an isomorphism
\begin{equation}\label{exp-8}
F^1H^2_\dR(\ol{U},D)\os{\cong}{\lra} F^1H^2_\dR(\ol{U})_0\os{\cong}{\longleftarrow}\goodxa.
\end{equation}
We denote by $\omega_{\ol{U},D}\in F^1H^2_\dR(\ol{U},D)$ the corresponding 
element of $\omega\in \goodxa$ via \eqref{exp-8}.
Fix a $\Gamma\in H_2(\ol{U^o},D,\Q)$ such that $\gamma=\partial(\Gamma)$.
Then,
by Thm. \ref{bregin-thm}, \eqref{beireg-3} and Prop. \ref{pairingVD4},
we have
\begin{equation}\label{exp-7}
\reg(\xi_\gamma)=
\left[\omega\longmapsto\int_\Gamma \omega_{\ol{U},D}
=\int_\Gamma \omega-\int_{\partial\Gamma}\cand_D(\omega)\right]
\in \Hom(\goodxa,\C)/\Image \Phi
\end{equation}
under the identification \eqref{exp-6}.

\subsection{$\bE(\ol{U^o},D;\Z)$ and $\bE(\ol{U^o},\Z)$}\label{exp-sect2}
Take a path $\gamma:[0,1]\to \ol{S^o}(\C)$, $t\mapsto\gamma_t$
such that $\gamma_t\in S^o(\C)$ for $t\ne0,1$.
Take a cycle $\varepsilon\in H_1(f^{-1}(\gamma_{t_0}),\Z)$ for some (fixed) $t_0\in [0,1]$.
Then it extends to a flat section 
$\varepsilon_{t}\in H_1(f^{-1}(\gamma_{t}),\Z)$ over $t\in [0,1]$
in a unique way.
We denote by $\Gamma(\varepsilon,\gamma)$ the fibration over the path $\gamma$ 
whose fiber is $\varepsilon_t$.

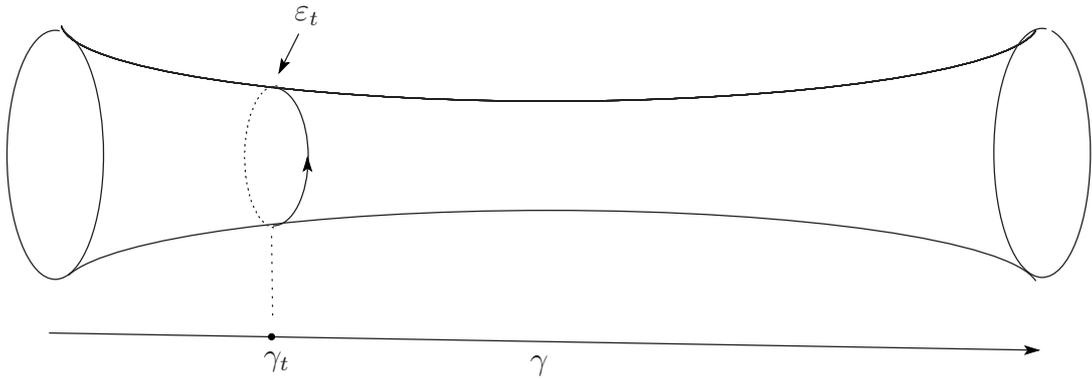
\unitlength 0.1in
\begin{picture}( 56.1000, 19.7500)(  5.1000,-34.7500)
%
{\color[named]{Black}{%
\special{pn 8}%
\special{ar 3320 1740 2530 410  0.1047219  3.0932053}%
}}%
%
{\color[named]{Black}{%
\special{pn 8}%
\special{ar 3320 1740 2530 410  0.1047219  3.0932053}%
}}%
%
{\color[named]{Black}{%
\special{pn 8}%
\special{ar 3320 1740 2530 410  0.1047219  3.0932053}%
}}%
%
{\color[named]{Black}{%
\special{pn 8}%
\special{ar 3320 1740 2530 410  0.1047219  3.0932053}%
}}%
%
{\color[named]{Black}{%
\special{pn 8}%
\special{ar 3320 1740 2530 410  0.1047219  3.0932053}%
}}%
%
{\color[named]{Black}{%
\special{pn 8}%
\special{ar 3320 1740 2530 410  0.1047219  3.0932053}%
}}%
%
{\color[named]{Black}{%
\special{pn 8}%
\special{ar 3320 1740 2530 410  0.1047219  3.0932053}%
}}%
%
{\color[named]{Black}{%
\special{pn 8}%
\special{ar 3320 1740 2530 410  0.1047219  3.0932053}%
}}%
%
{\color[named]{Black}{%
\special{pn 8}%
\special{ar 3320 1740 2530 410  0.1047219  3.0932053}%
}}%
%
{\color[named]{Black}{%
\special{pn 8}%
\special{ar 3320 1740 2530 410  0.1047219  3.0932053}%
}}%
%
{\color[named]{Black}{%
\special{pn 8}%
\special{ar 3320 3130 2530 410  3.3177741  6.1784634}%
}}%
%
{\color[named]{Black}{%
\special{pn 8}%
\special{ar 5870 2420 250 650  4.9194617  4.7935114}%
}}%
%
{\color[named]{Black}{%
\special{pn 8}%
\special{ar 760 2430 250 650  4.9194617  4.7935114}%
}}%
%
{\color[named]{Black}{%
\special{pn 8}%
\special{ar 1890 2440 180 360  4.6583875  1.5707963}%
}}%
%
{\color[named]{Black}{%
\special{pn 8}%
\special{pa 1890 2810}%
\special{pa 1884 2808}%
\special{fp}%
\special{pa 1862 2800}%
\special{pa 1856 2796}%
\special{fp}%
\special{pa 1838 2782}%
\special{pa 1834 2776}%
\special{fp}%
\special{pa 1818 2758}%
\special{pa 1814 2752}%
\special{fp}%
\special{pa 1802 2732}%
\special{pa 1800 2726}%
\special{fp}%
\special{pa 1788 2704}%
\special{pa 1786 2698}%
\special{fp}%
\special{pa 1776 2674}%
\special{pa 1774 2668}%
\special{fp}%
\special{pa 1766 2644}%
\special{pa 1766 2638}%
\special{fp}%
\special{pa 1760 2614}%
\special{pa 1758 2606}%
\special{fp}%
\special{pa 1752 2582}%
\special{pa 1752 2576}%
\special{fp}%
\special{pa 1748 2550}%
\special{pa 1746 2542}%
\special{fp}%
\special{pa 1744 2518}%
\special{pa 1744 2510}%
\special{fp}%
\special{pa 1742 2484}%
\special{pa 1742 2476}%
\special{fp}%
\special{pa 1740 2450}%
\special{pa 1740 2444}%
\special{fp}%
\special{pa 1740 2416}%
\special{pa 1742 2410}%
\special{fp}%
\special{pa 1742 2384}%
\special{pa 1742 2376}%
\special{fp}%
\special{pa 1746 2350}%
\special{pa 1746 2344}%
\special{fp}%
\special{pa 1750 2318}%
\special{pa 1750 2310}%
\special{fp}%
\special{pa 1756 2286}%
\special{pa 1756 2278}%
\special{fp}%
\special{pa 1762 2254}%
\special{pa 1764 2248}%
\special{fp}%
\special{pa 1770 2224}%
\special{pa 1772 2218}%
\special{fp}%
\special{pa 1782 2194}%
\special{pa 1784 2188}%
\special{fp}%
\special{pa 1794 2164}%
\special{pa 1796 2158}%
\special{fp}%
\special{pa 1808 2138}%
\special{pa 1812 2132}%
\special{fp}%
\special{pa 1826 2112}%
\special{pa 1830 2108}%
\special{fp}%
\special{pa 1846 2092}%
\special{pa 1852 2088}%
\special{fp}%
\special{pa 1872 2076}%
\special{pa 1878 2074}%
\special{fp}%
\special{pa 1902 2070}%
\special{pa 1910 2072}%
\special{fp}%
}}%
%
{\color[named]{Black}{%
\special{pn 8}%
\special{pa 2020 1800}%
\special{pa 1920 2000}%
\special{fp}%
\special{sh 1}%
\special{pa 1920 2000}%
\special{pa 1968 1950}%
\special{pa 1944 1952}%
\special{pa 1932 1932}%
\special{pa 1920 2000}%
\special{fp}%
}}%
%
{\color[named]{Black}{%
\special{pn 8}%
\special{pa 730 3360}%
\special{pa 5850 3450}%
\special{fp}%
\special{sh 1}%
\special{pa 5850 3450}%
\special{pa 5784 3430}%
\special{pa 5798 3450}%
\special{pa 5784 3470}%
\special{pa 5850 3450}%
\special{fp}%
}}%
\put(32.2000,-35.8000){\makebox(0,0)[lb]{$\gamma$}}%
\put(20.0000,-17.5000){\makebox(0,0)[lb]{$\varepsilon_t$}}%
\put(18.4000,-35.6000){\makebox(0,0)[lb]{$\gamma_t$}}%
%
{\color[named]{Black}{%
\special{pn 8}%
\special{pa 1880 2856}%
\special{pa 1886 3300}%
\special{dt 0.045}%
}}%
%
{\color[named]{Black}{%
\special{pn 4}%
\special{sh 1}%
\special{ar 1880 3380 16 16 0  6.28318530717959E+0000}%
\special{sh 1}%
}}%
%
{\color[named]{Black}{%
\special{pn 8}%
\special{pa 2066 2466}%
\special{pa 2066 2446}%
\special{fp}%
\special{sh 1}%
\special{pa 2066 2446}%
\special{pa 2046 2512}%
\special{pa 2066 2498}%
\special{pa 2086 2512}%
\special{pa 2066 2446}%
\special{fp}%
}}%
\end{picture}%

\vspace{0.5cm}

Then \[
\Gamma(\varepsilon,\gamma)\in H_2(\ol{U^o},f^{-1}(\gamma_0)\cup f^{-1}(\gamma_1);\Z),
\]
\[
\mbox{with }\quad
\partial(\Gamma(\varepsilon,\gamma))=\varepsilon_1-\varepsilon_0
\in H_1(f^{-1}(\gamma_0)\cup f^{-1}(\gamma_1),\Z).
\]
Define
$
\bE(\ol{U^o},D;\Z)\subset H_2(\ol{U^o},D;\Z)
$
the subgroup generated by $\Gamma(\varepsilon,\gamma)$'s where $\gamma$ and $\varepsilon$
run over as above such that $\gamma_0,\gamma_1\in T_m=\ol{S^o}-S^o$.
Define $\bE(\ol{U^o},\Z)$ by an exact sequence
\[
0\lra \bE(\ol{U^o},\Z)\lra \bE(\ol{U^o},D;\Z)\os{\partial}{\lra} H_1(D,\Z).
\]
Write $\bE(\ol{U^o},D;\Q):=\bE(\ol{U^o},D;\Z)\ot\Q$ etc.
\begin{prop}\label{Ev}We have
\begin{equation}\label{exp-4}
\bE(\ol{U^o},D;\Q)= H_2(\ol{U^o},D;\Q),
\end{equation}
\begin{equation}\label{exp-5}
\bE(\ol{U^o},\Q)\os{\cong}{\lra} H_2(\ol{U^o},\Q)/H_2(D,\Q)\cong (H^2(\ol{U^o})_0)^*.
\end{equation}
Hence we have 
\[
\xymatrix{
0\ar[r]&\bE(\ol{U^o};\Q)\ar[r]\ar[d]^\cong&\bE(\ol{U^o},D;\Q)\ar[r]^\partial
\ar@{=}[d]&
H_1(D,\Q)\ar@{=}[d]\ar[r]&0\\
0\ar[r]&H_2(\ol{U^o},\Q)/H_2(D,\Q)\ar[r]&H_2(\ol{U^o},D;\Q)\ar[r]^\partial&H_1(D,\Q)\ar[r]&0.
}
\]
\end{prop}

\begin{lem}\label{Ev-lem}
The sequence
\begin{equation}\label{Ev-lem-2}
H_2(U^o,\Q)\lra H_2(\ol{U^o},D;\Q)\os{\partial}{\lra} H_1(D,\Q)\lra 0
\end{equation}
is exact. 
\end{lem}
\begin{pf}
The surjectivity of $\partial$ is proven in the same way as Lem. \ref{exp-lem0}.
We only show
\begin{equation}\label{Ev-lem-1}
\Image(H_2(U^o,\Q)\lra H_2(\ol{U^o},D;\Q))
=\Image(H_2(\ol{U^o},\Q)\lra H_2(\ol{U^o},D;\Q)).
\end{equation}
Consider a diagram 
\[
\xymatrix{
&H_2(D,\Q)\ar[d]^a\\
H_2(U^o,\Q)\ar[r]& H_2(\ol{U^o},\Q)\ar[d]\ar[r]^b& H^2(D,\Q)\ar[r]
&H_1(U^o,\Q)\ar[r]^c&H_1(\ol{U^o},\Q)\\
&H_2(\ol{U^o},D;\Q)\\
}\]
with exact row and column.
Hence it is enough to show $\Image(ba)=\Image(b)$ or equivalently $\dim\Coker(ba)=\dim\Coker(b)
(=\dim\ker(c))$.
Since $ba$ is given by the intersection pairing, a direct calculation shows that
$\dim\Coker(ba)=\dim H_0(T_m)$.
On the other hand, \[\ker(H_1(U^o,\Q)\os{c}{\to} H_1(\ol{U^o},\Q))\cong
\ker(H_1(S^o,\Q)\to H_1(\ol{S^o},\Q))\cong H_0(T_m)\] (cf. \eqref{exp-lem0-1}), 
so we are done.
\end{pf}
\noindent{\it Proof of Prop.\ref{Ev}.}
Let $\cL$ be the local system on $S^o(\C)$ whose fiber is $H_1(f^{-1}(t),\Q)$.
Then the image of $H_2(U,\Q)$ in $H_2(\ol{U^o},D;\Q)$ coincides with that of
$H_1(S^o,\cL)$. The homology group $H_1(S^o,\cL)$ is generated by
$\Gamma(\varepsilon,\gamma)$'s such that $\gamma_0=\gamma_1\in S^o$ and $\varepsilon_0=
\varepsilon_1$. Take a path $\delta$ such that $\delta_0\in T_m=\ol{S^o}-S^o$ 
and $\delta_1=\gamma_0=\gamma_1$.
Put $\wt{\gamma}=\delta\cdot \gamma\cdot\delta^{-1}$.
Then the image of $\Gamma(\varepsilon,\gamma)$ in $H_2(\ol{U^o},D;\Q)$ coincides with
$\Gamma(\varepsilon,\wt{\gamma})$, and this is an element of $\bE(\ol{U^o},D;\Q)$.
There remains to show the surjectivity of $\bE(\ol{U^o},D;\Q)\to H_1(D,\Q)$
(this gives an alternative proof of Lem.\ref{exp-lem0}).
To do this, it is enough to show that for each $p\in T_m$, there is a path $\nu$
such that $\nu_0=p$ and $\nu_1\in T_m$,
and there is a cycle $\alpha_t\in H_1(f^{-1}(\nu_t),\Q)$ such that
$\alpha_0\ne0$ and $\alpha_1=0$.
Since
\[
\partial(\Gamma(\alpha,\nu))=(\cdots,0,\alpha_0,0,\cdots)\in H_1(D,\Q)=\bigoplus_{x\in T_m}
H_1(f^{-1}(x),\Q)
\]
this implies the surjectivity of $\bE(\ol{U^o},D;\Q)\to H_1(D,\Q)$.
Fix paths $\nu'$ and $\nu^{\prime\prime}$ such that $\nu'_0=p$,
$\nu'_1=\nu^{\prime\prime}_1=q\in S^o$ and $\nu^{\prime\prime}_0\in T_m$.
Fix $\alpha'\in H_1(f^{-1}(q),\Q)$ such that $\alpha'$ goes to a nonzero cycle
as $q\to p$,
and $\alpha^{\prime\prime}\in H_1(f^{-1}(q),\Q)$ such that $\alpha^{\prime\prime}$
goes to zero
as $q\to \nu^{\prime\prime}_0$.
Since $H_1(f^{-1}(q),\Q)$ is an irreducible $\Q[\pi_1(S^o,q)]$-module,
there is $g\in \Q[\pi_1(S^o,q)]$ such that $g(\alpha')=\alpha^{\prime\prime}$.
Joining $\Gamma(\alpha',\nu')$, $\Gamma(\alpha',g)$ and 
$\Gamma(\alpha^{\prime\prime},\nu^{\prime\prime})$, we obtain 
$\Gamma(\alpha,\nu)$ as desired.
This completes the proof of Prop.\ref{Ev}.

\subsection{A formula of Beilinson regulator on $H^3_{\cM}(X,\Q(2))$}\label{exp-sect3}
We summarize all of the results in \S \ref{exp-sect1} and \S \ref{exp-sect2} together with
Thm.\ref{bregin-thm} in the following theorem.
\begin{thm}\label{ExtMHS}
Let the notations be as in \S \ref{deRh-sect-1}, \S \ref{good-sect} and \S \ref{exp-sect2}.
Suppose $K=\C$ and $D\ne\emptyset$.
Then we have
\[
\Ext^1(\Q,\NF(X)^\perp\ot\Q(2))\cong\Coker
\left[\Phi:
\bE(\ol{U^o},\Q)\lra \Hom(\goodxa,\C)
\right]
\]
where 
\[
\Phi(\Delta)=\left[\omega\longmapsto\int_\Delta \omega\right],\quad \omega\in \goodxa.
\]
For $\gamma\in H_1(D,\Q)$, fix $\Gamma\in \bE(\ol{U^o},D;\Q)$
with $\partial(\Gamma)=\gamma$.
Then 
\[
\reg(\xi_\gamma)=\pm
\left[\omega\longmapsto\int_\Gamma \omega_{\ol{U^o},D}
=\int_\Gamma \omega-\int_{\gamma}\cand_D(\omega)\right]\in \Hom(\goodxa,\C)/\Image\Phi
\]
(see Prop. \ref{pairingVD4} and \eqref{pairingVD3} for ``$\cand_D$").
Moreover if $\gamma\in H_1(D_0,\Q)$ with $D_0$ a union of irreducible multiplicative fibers,
then $\cand_{D_0}(\omega)=0$ by Prop. \ref{hodge-mot}
so that we have simply
\[
\reg(\xi_\gamma)=\pm
\left[\omega\longmapsto\int_\Gamma \omega\right].
\]
\end{thm}
The point is that ``$\omega\in \goodxa$" is an algebraic 2-form.
This makes it easier to compute the regulator.
To carry out the computation practically, we need the following data.
\begin{itemize}
\item[(a)]
A basis of $\goodxa$,
\item[(b)]
A basis of $\ol{H_2(X,D;\Q)}$ (see \eqref{main-diagram}),
\item[(c)]
Computation of the extra term $\cand_D$ (however see Rem. \ref{hodge-mot-rem}).
\end{itemize}
One can obtain (a) by a direct computation of $H^1_\dR(C,\cH_e)$
and by explicit formula of Gauss-Manin connection (Appendix).
See \S \ref{good-expm-sect} for an example of the computation. 
A basis (b) can be constructed from 
$H_2(\ol{U^o},D;\Q)\cong\bE(\ol{U^o},D;\Q)$.
It is not hard to obtain
a basis of $\bE(\ol{U^o},D;\Q)/\bE(\ol{U^o},\Q)$.
To obtain a basis of $H_2(X,\Q)/\NF^B(X)$
we assume that the precision of the values of integrations 
can be raised as many as one likes.
Then, by using the basis of $\goodxa$ together with the fact that
there is an embedding
$\Image[\bE(\ol{U^o},\Q)\to H_2(X,\Q)/\NF^B(X)]\hra \Hom(\goodxa,\C)$,
one can prove the linear independence of given cycles in $\bE(\ol{U^o},\Q)$
if they were linear independent.
Hence one can eventually obtain a basis of $\Image\bE(\ol{U^o},\Q)$.

\medskip

We shall apply the above method to an example in the next section.
\section{Example : $3y^2+x^3+(3x+4t^l)^2$}\label{Example-sect}
Let $l\geq1$ be an integer. We consider a minimal elliptic surface
\[
f:X\lra \P^1,\quad f^{-1}(t):3y^2+x^3+(3x+4t^l)^2=0
\]
defined over $\Q$. There is the section $e:\P^1\to X$ of ``$\infty$".
Write $X_\C:=X\times_\Q\C$.

\medskip

The purpose of this section is to compute the real regulator
\[
\reg_\R:H^3_{\mathscr M}(X,\Q(2))\lra
\mathrm{Ext}^1_{\R\mbox{-}\MHS}(\R,H^2(X)_\ind\ot\R(1))^{F_\infty=1},
\]
where $H^2(X)_\ind:=H^2(X_\C,\Q(1))/\NS(X_\C)$, especially for an element
\[
\xi_{D_1}=\left[\frac{y-(x+4)}{y+(x+4)},D_1\right]\in H^3_{\mathscr M}(X,\Q(2))\quad
(D_1:=f^{-1}(1))
\]
arising from a split multiplicative fiber $D_1$ of type $\mathrm{I}_1$.
We note that if $(l,6)=1$ then 
$\xi_{D_1}$ is an ``integral" element, in the sense that it comes from the motivic
cohomology of a proper flat regular model of $X$ over $\Z$ (see \cite{scholl} 1.1.6 for 
``unconditional" definition of integral elements).
\subsection{Basic data of $X$}
The following is easy to show (the proof is left to the reader).
\begin{itemize}
\item The Hodge numbers are as follows:
\[
h^{10}=h^{01}=0,\quad
h^{20}=h^{02}=
\lfloor\frac{l-1}{3}\rfloor,\quad
h^{11}=10(1+h^{20}).\]
In particular, $H^2_B(X)_\ind:=H^2_B(X(\C),\Q)/\NS(X_\C)\ot\Q\ne0$ if and only if $l\geq 4$.
(If $1\leq l\leq 3$, then $X$ is a rational surface.)  
\item
There are $(l+1)$-multiplicative fibers: 
\[D_0:=f^{-1}(0) \mbox{(=type }\mathrm{I}_{3l}),\quad 
D_i:=f^{-1}(\zeta_l^{i-1}) \mbox{(=type }\mathrm{I}_1),\]
where $\zeta_l=\exp(2\pi i/l)$ and $1\leq i \leq l$. 
Moreover $D_1=f^{-1}(1)$ is the unique split multiplicative fiber.
\item
If $3|l$, then there is no additive fiber.
If $(3,l)=1$ then $E=f^{-1}(\infty)$ is 
the unique additive fiber (type IV* if $l\equiv 1$ mod 3, and type IV
if $l\equiv2$ mod 3).
In particular, $E\ne\emptyset$ if and only if $(l,3)=1$.
\item $\NF(X_\C)\ot\Q=\NS(X_\C)\ot\Q$ if and only if $l$ is odd
(\cite{stiller2} Example 4).
\item
There is an automorphism $\sigma:X\to X$ given by $\sigma(x,y,t)=(x,y,\zeta_lt)$.
\end{itemize}


Hereafter we assume $(l,6)=1$.
Then 
\[
\NF_B(X_\C)^\perp
=\NS(X_\C)^\perp
\os{\cong}{\lra}
H^2_B(X)_\ind:=H^2_B(X(\C),\Q(1))/\NS(X_\C)\ot\Q.
\]
Since $U=U^o$ in this case one has
\begin{equation}\label{Example-7}
\dim \bE(\ol{U},\Q)=\dim \NF(X)^\perp=l-1\mbox{ (by \eqref{exp-5} and Prop.\ref{VDR} (3))},
\end{equation}
\begin{equation}\label{Example-8}
\dim \bE(\ol{U},D;\Q)=(l-1)+\dim H_1(D,\Q)=2l-1.
\end{equation}

\subsection{Good algebraic 2-forms $\goodxa$ and $\goodxb$}\label{good-expm-sect}
Suppose $(l,6)=0$.
We use the same notations in \S \ref{deRh-sect-1},
$D=D_0+\cdots+D_l$, 
$U=X-(D+E)$, $\ol{U}=X-E$, $T_m=\{0,\zeta_l^i\}$, $T_a=\{\infty\}$ and $T=T_m+T_a$. 
Let $H^2_\dR(\ol{U})_0:=\ker[H^2_\dR(\ol{U})\to H^2_\dR(D)]$.
This is isomorphic to $\NF_\dR(X)^\perp$ by Prop. \ref{VDR} (3).

\medskip

As is easily seen, one has
\begin{equation}
\goodxb=\langle t^{i-1}\frac{dtdx}{y}~|~1\leq i\leq \lfloor\frac{l-1}{3}\rfloor
\rangle_\Q\os{\cong}{\lra}F^2 H^2_\dR(\ol{U})_0.
\end{equation}
Let us compute $\goodxa$.
Since $H^1_\zar(\P^1,\cH_e^{1,0})=0$, $\goodua$ is the image of the composition
\begin{align*}
\vg(\P^1,\Omega^1_{\P^1}(\log T)\ot\cH_e)&\lra 
\vg(\bA^1,\Omega^1_{\P^1}(\log T)\ot\cH_e)/\nabla(\cH^{1,0}_e))\\
&\os{\cong}{\longleftarrow}
\vg(\bA^1,\Omega^1_{\P^1}(\log T)\ot\cH^{1,0}_e)\\
&=\vg(\bA^1,\Omega^2_{\bA^1}(\log T_m))
\end{align*}
where $\bA^1=\P^1-\{\infty\}$.
Using the basis $\canh$ and $\can$ in Appendix \eqref{A1-1}, \eqref{A1-2}, 
one easily sees that $\vg(\P^1,\Omega^1_{\P^1}(\log T)\ot\cH_e)$
is generated by the following elements.
\[
\frac{t^idt}{t(t^l-1)}\ot\canh~(0\leq i\leq l+\lfloor\frac{l-1}{3}\rfloor),\quad
\frac{t^jdt}{t(t^l-1)}\ot\can~(0\leq j\leq l-1-\lfloor\frac{l-1}{3}\rfloor).
\]
By using Thm. \ref{A1} \eqref{A1-3} and \eqref{A1-4}, we can compute
their image into 
$\vg(\bA^1,\Omega^2_{\bA^1}(\log T_m))$ directly.
We thus have
\begin{equation}
\goodua=\langle \frac{t^idt}{t(t^l-1)}\frac{dx}{y}
~|~0\leq i\leq 2l-1-\lfloor\frac{l-1}{3}\rfloor\rangle_\Q.
\end{equation}
Therefore
\begin{equation}
\goodxa=\langle t^{i-1}\frac{dtdx}{y}~|~1\leq i\leq l-1-\lfloor\frac{l-1}{3}\rfloor
\rangle_\Q
\end{equation}
By Prop. \ref{hodge-mot}, the extra term for $D_k$ $(1\leq k\leq l)$ vanishes:
\begin{equation}
\cand_{D_k} \left(t^{i-1}\frac{dtdx}{y}\right)=0~\mbox{ for }
1\leq i\leq l-1-\lfloor\frac{l-1}{3}\rfloor,~1\leq k\leq l.
\end{equation}
\begin{exmp}
One can show that 
\[
(l+3i)(2l+3i)t^{i-1+l}\frac{dtdx}{y}
\equiv9i^2t^{i-1}\frac{dtdx}{y}
\]
in $H^2_\dR(\ol{U})$, however
\[
\int_\epsilon\cand_{D_1}\left((l+3i)(2l+3i)t^{i-1+l}\frac{dtdx}{y}
-9i^2t^{i-1}\frac{dtdx}{y}
\right)=\pm9l\ne0
\]
by using Prop. \ref{extraA1}, \eqref{extraA3}, 
where $\epsilon\in H_1(D_1,\Z)\cong\Z$ is the generator.
In particular $\cand_{D_1}(t^{i-1+l}\frac{dtdx}{y})\ne0$.
\end{exmp}


\subsection{Cycles $\Delta$ and $\Gamma$}
Let $\delta_0$ (resp. $\delta_1$) be the homology cycle in $H_1(f^{-1}(t),\Z)$
which vanishes as $t\to 0$ (resp. $t\to 1$).
Define $\Delta$ and $\Gamma$ to be fibrations
over the segment $[0,1]\subset \P^1(\C)$ whose fibers are the vanishing cycles
$\delta_1$ and $\delta_0$ respectively.
\[
\Delta\in H_2(\ol{U},D_0;\Z),\quad
\Gamma\in H_2(\ol{U},D_1;\Z).
\]
The boundary $\partial\Delta$ (resp. $\partial\Gamma$) is a generator
of the homology group $H_1(D_0,\Z)\cong \Z$ (resp. $H_1(D_1,\Z)\cong \Z$).

\vspace{0.5cm}

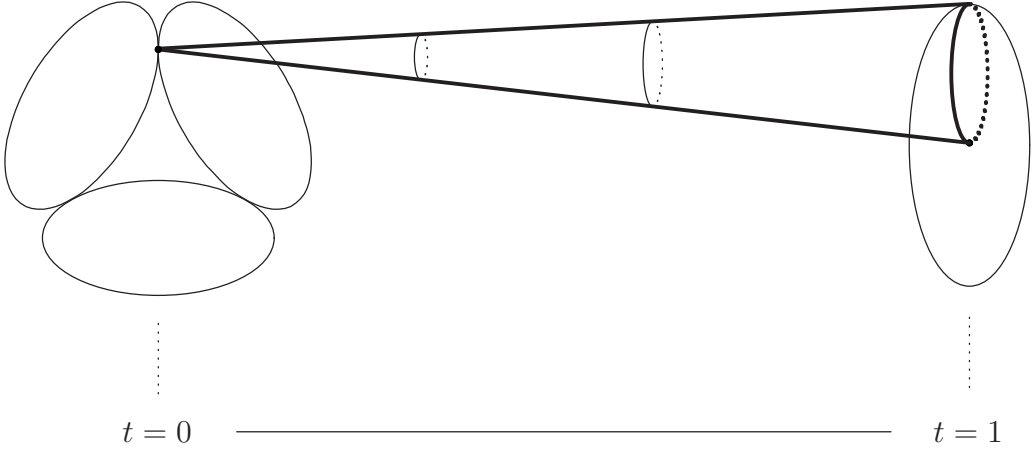
\begin{center}
\unitlength 0.1in
\begin{picture}( 53.0500, 53.6500)(  9.0200,-62.0000)
%
{\color[named]{Black}{%
\special{pn 8}%
\special{pa 5896 2602}%
\special{pa 5896 3002}%
\special{dt 0.045}%
\special{pa 1696 3002}%
\special{pa 1696 2602}%
\special{dt 0.045}%
}}%
%
{\color[named]{Black}{%
\special{pn 8}%
\special{ar 5896 4704 312 736  0.0000000  6.2831853}%
}}%
%
{\color[named]{Black}{%
\special{pn 4}%
\special{sh 1}%
\special{ar 5896 4692 8 8 0  6.28318530717959E+0000}%
\special{sh 1}%
\special{ar 5896 4692 8 8 0  6.28318530717959E+0000}%
}}%
%
{\color[named]{Black}{%
\special{pn 8}%
\special{pa 5896 5602}%
\special{pa 5896 6002}%
\special{dt 0.045}%
\special{pa 1696 6002}%
\special{pa 1696 5602}%
\special{dt 0.045}%
}}%
\put(16.9000,-62.0000){\makebox(0,0){$t=0$}}%
\put(58.9000,-62.0000){\makebox(0,0){$t=1$}}%
%
{\color[named]{Black}{%
\special{pn 8}%
\special{pa 5490 6200}%
\special{pa 2102 6200}%
\special{fp}%
}}%
%
{\color[named]{Black}{%
\special{pn 8}%
\special{ar 1696 5188 600 300  0.0000000  6.2831853}%
}}%
%
{\color[named]{Black}{%
\special{pn 8}%
\special{pa 1598 3978}%
\special{pa 1624 3998}%
\special{pa 1646 4020}%
\special{pa 1664 4048}%
\special{pa 1676 4078}%
\special{pa 1686 4108}%
\special{pa 1692 4140}%
\special{pa 1696 4204}%
\special{pa 1692 4268}%
\special{pa 1688 4298}%
\special{pa 1682 4330}%
\special{pa 1674 4362}%
\special{pa 1666 4392}%
\special{pa 1658 4424}%
\special{pa 1638 4484}%
\special{pa 1626 4514}%
\special{pa 1612 4542}%
\special{pa 1598 4572}%
\special{pa 1570 4628}%
\special{pa 1554 4656}%
\special{pa 1536 4684}%
\special{pa 1520 4710}%
\special{pa 1502 4738}%
\special{pa 1482 4764}%
\special{pa 1464 4788}%
\special{pa 1442 4814}%
\special{pa 1422 4836}%
\special{pa 1378 4884}%
\special{pa 1330 4928}%
\special{pa 1306 4948}%
\special{pa 1280 4966}%
\special{pa 1224 4998}%
\special{pa 1196 5012}%
\special{pa 1166 5024}%
\special{pa 1136 5032}%
\special{pa 1104 5038}%
\special{pa 1072 5038}%
\special{pa 1040 5034}%
\special{pa 1010 5024}%
\special{pa 982 5006}%
\special{pa 958 4984}%
\special{pa 940 4958}%
\special{pa 926 4930}%
\special{pa 916 4900}%
\special{pa 908 4868}%
\special{pa 904 4836}%
\special{pa 902 4804}%
\special{pa 902 4772}%
\special{pa 904 4740}%
\special{pa 912 4676}%
\special{pa 920 4646}%
\special{pa 928 4614}%
\special{pa 936 4584}%
\special{pa 956 4524}%
\special{pa 980 4464}%
\special{pa 994 4434}%
\special{pa 1022 4378}%
\special{pa 1054 4322}%
\special{pa 1072 4294}%
\special{pa 1108 4242}%
\special{pa 1128 4218}%
\special{pa 1148 4192}%
\special{pa 1168 4168}%
\special{pa 1212 4120}%
\special{pa 1234 4098}%
\special{pa 1258 4076}%
\special{pa 1282 4056}%
\special{pa 1308 4038}%
\special{pa 1336 4020}%
\special{pa 1362 4004}%
\special{pa 1390 3988}%
\special{pa 1420 3976}%
\special{pa 1450 3966}%
\special{pa 1482 3958}%
\special{pa 1514 3956}%
\special{pa 1546 3958}%
\special{pa 1576 3968}%
\special{pa 1598 3978}%
\special{fp}%
}}%
%
{\color[named]{Black}{%
\special{pn 8}%
\special{pa 1792 3978}%
\special{pa 1766 3998}%
\special{pa 1744 4020}%
\special{pa 1728 4048}%
\special{pa 1714 4078}%
\special{pa 1706 4108}%
\special{pa 1700 4140}%
\special{pa 1696 4204}%
\special{pa 1696 4236}%
\special{pa 1698 4268}%
\special{pa 1702 4298}%
\special{pa 1708 4330}%
\special{pa 1716 4362}%
\special{pa 1724 4392}%
\special{pa 1734 4424}%
\special{pa 1754 4484}%
\special{pa 1766 4514}%
\special{pa 1778 4542}%
\special{pa 1792 4572}%
\special{pa 1806 4600}%
\special{pa 1854 4684}%
\special{pa 1872 4710}%
\special{pa 1890 4738}%
\special{pa 1908 4764}%
\special{pa 1928 4788}%
\special{pa 1948 4814}%
\special{pa 1970 4836}%
\special{pa 2014 4884}%
\special{pa 2036 4906}%
\special{pa 2060 4928}%
\special{pa 2086 4948}%
\special{pa 2112 4966}%
\special{pa 2138 4982}%
\special{pa 2166 4998}%
\special{pa 2196 5012}%
\special{pa 2226 5024}%
\special{pa 2256 5032}%
\special{pa 2288 5038}%
\special{pa 2320 5038}%
\special{pa 2352 5034}%
\special{pa 2382 5024}%
\special{pa 2408 5006}%
\special{pa 2432 4984}%
\special{pa 2452 4958}%
\special{pa 2466 4930}%
\special{pa 2476 4900}%
\special{pa 2482 4868}%
\special{pa 2486 4836}%
\special{pa 2488 4804}%
\special{pa 2488 4772}%
\special{pa 2484 4708}%
\special{pa 2478 4676}%
\special{pa 2472 4646}%
\special{pa 2464 4614}%
\special{pa 2454 4584}%
\special{pa 2446 4554}%
\special{pa 2398 4434}%
\special{pa 2384 4406}%
\special{pa 2368 4378}%
\special{pa 2354 4350}%
\special{pa 2338 4322}%
\special{pa 2320 4294}%
\special{pa 2302 4268}%
\special{pa 2282 4242}%
\special{pa 2264 4218}%
\special{pa 2244 4192}%
\special{pa 2224 4168}%
\special{pa 2180 4120}%
\special{pa 2158 4098}%
\special{pa 2134 4076}%
\special{pa 2108 4056}%
\special{pa 2056 4020}%
\special{pa 2000 3988}%
\special{pa 1970 3976}%
\special{pa 1940 3966}%
\special{pa 1908 3958}%
\special{pa 1876 3956}%
\special{pa 1844 3958}%
\special{pa 1814 3968}%
\special{pa 1792 3978}%
\special{fp}%
}}%
%
{\color[named]{Black}{%
\special{pn 8}%
\special{ar 5896 1704 312 736  0.0000000  6.2831853}%
}}%
%
{\color[named]{Black}{%
\special{pn 4}%
\special{sh 1}%
\special{ar 5896 1692 8 8 0  6.28318530717959E+0000}%
\special{sh 1}%
\special{ar 5896 1692 8 8 0  6.28318530717959E+0000}%
}}%
%
{\color[named]{Black}{%
\special{pn 8}%
\special{pa 5896 2602}%
\special{pa 5896 3002}%
\special{dt 0.045}%
\special{pa 1696 3002}%
\special{pa 1696 2602}%
\special{dt 0.045}%
}}%
%
{\color[named]{Black}{%
\special{pn 8}%
\special{ar 1696 2188 600 300  0.0000000  6.2831853}%
}}%
%
{\color[named]{Black}{%
\special{pn 8}%
\special{pa 1598 978}%
\special{pa 1624 998}%
\special{pa 1646 1020}%
\special{pa 1664 1048}%
\special{pa 1676 1078}%
\special{pa 1686 1108}%
\special{pa 1692 1140}%
\special{pa 1696 1204}%
\special{pa 1692 1268}%
\special{pa 1688 1298}%
\special{pa 1682 1330}%
\special{pa 1674 1362}%
\special{pa 1666 1392}%
\special{pa 1658 1424}%
\special{pa 1638 1484}%
\special{pa 1626 1514}%
\special{pa 1612 1542}%
\special{pa 1598 1572}%
\special{pa 1570 1628}%
\special{pa 1554 1656}%
\special{pa 1536 1684}%
\special{pa 1520 1710}%
\special{pa 1502 1738}%
\special{pa 1482 1764}%
\special{pa 1464 1788}%
\special{pa 1442 1814}%
\special{pa 1422 1836}%
\special{pa 1378 1884}%
\special{pa 1330 1928}%
\special{pa 1306 1948}%
\special{pa 1280 1966}%
\special{pa 1224 1998}%
\special{pa 1196 2012}%
\special{pa 1166 2024}%
\special{pa 1136 2032}%
\special{pa 1104 2038}%
\special{pa 1072 2038}%
\special{pa 1040 2034}%
\special{pa 1010 2024}%
\special{pa 982 2006}%
\special{pa 958 1984}%
\special{pa 940 1958}%
\special{pa 926 1930}%
\special{pa 916 1900}%
\special{pa 908 1868}%
\special{pa 904 1836}%
\special{pa 902 1804}%
\special{pa 902 1772}%
\special{pa 904 1740}%
\special{pa 912 1676}%
\special{pa 920 1646}%
\special{pa 928 1614}%
\special{pa 936 1584}%
\special{pa 956 1524}%
\special{pa 980 1464}%
\special{pa 994 1434}%
\special{pa 1022 1378}%
\special{pa 1054 1322}%
\special{pa 1072 1294}%
\special{pa 1108 1242}%
\special{pa 1128 1218}%
\special{pa 1148 1192}%
\special{pa 1168 1168}%
\special{pa 1212 1120}%
\special{pa 1234 1098}%
\special{pa 1258 1076}%
\special{pa 1282 1056}%
\special{pa 1308 1038}%
\special{pa 1336 1020}%
\special{pa 1362 1004}%
\special{pa 1390 988}%
\special{pa 1420 976}%
\special{pa 1450 966}%
\special{pa 1482 958}%
\special{pa 1514 956}%
\special{pa 1546 958}%
\special{pa 1576 968}%
\special{pa 1598 978}%
\special{fp}%
}}%
%
{\color[named]{Black}{%
\special{pn 8}%
\special{pa 1792 978}%
\special{pa 1766 998}%
\special{pa 1744 1020}%
\special{pa 1728 1048}%
\special{pa 1714 1078}%
\special{pa 1706 1108}%
\special{pa 1700 1140}%
\special{pa 1696 1204}%
\special{pa 1696 1236}%
\special{pa 1698 1268}%
\special{pa 1702 1298}%
\special{pa 1708 1330}%
\special{pa 1716 1362}%
\special{pa 1724 1392}%
\special{pa 1734 1424}%
\special{pa 1754 1484}%
\special{pa 1766 1514}%
\special{pa 1778 1542}%
\special{pa 1792 1572}%
\special{pa 1806 1600}%
\special{pa 1854 1684}%
\special{pa 1872 1710}%
\special{pa 1890 1738}%
\special{pa 1908 1764}%
\special{pa 1928 1788}%
\special{pa 1948 1814}%
\special{pa 1970 1836}%
\special{pa 2014 1884}%
\special{pa 2036 1906}%
\special{pa 2060 1928}%
\special{pa 2086 1948}%
\special{pa 2112 1966}%
\special{pa 2138 1982}%
\special{pa 2166 1998}%
\special{pa 2196 2012}%
\special{pa 2226 2024}%
\special{pa 2256 2032}%
\special{pa 2288 2038}%
\special{pa 2320 2038}%
\special{pa 2352 2034}%
\special{pa 2382 2024}%
\special{pa 2408 2006}%
\special{pa 2432 1984}%
\special{pa 2452 1958}%
\special{pa 2466 1930}%
\special{pa 2476 1900}%
\special{pa 2482 1868}%
\special{pa 2486 1836}%
\special{pa 2488 1804}%
\special{pa 2488 1772}%
\special{pa 2484 1708}%
\special{pa 2478 1676}%
\special{pa 2472 1646}%
\special{pa 2464 1614}%
\special{pa 2454 1584}%
\special{pa 2446 1554}%
\special{pa 2398 1434}%
\special{pa 2384 1406}%
\special{pa 2368 1378}%
\special{pa 2354 1350}%
\special{pa 2338 1322}%
\special{pa 2320 1294}%
\special{pa 2302 1268}%
\special{pa 2282 1242}%
\special{pa 2264 1218}%
\special{pa 2244 1192}%
\special{pa 2224 1168}%
\special{pa 2180 1120}%
\special{pa 2158 1098}%
\special{pa 2134 1076}%
\special{pa 2108 1056}%
\special{pa 2056 1020}%
\special{pa 2000 988}%
\special{pa 1970 976}%
\special{pa 1940 966}%
\special{pa 1908 958}%
\special{pa 1876 956}%
\special{pa 1844 958}%
\special{pa 1814 968}%
\special{pa 1792 978}%
\special{fp}%
}}%
%
{\color[named]{Black}{%
\special{pn 20}%
\special{ar 1696 1730 520 520  0.0000000  6.2831853}%
}}%
%
{\color[named]{Black}{%
\special{pn 20}%
\special{pa 1694 1208}%
\special{pa 5896 1692}%
\special{fp}%
}}%
%
{\color[named]{Black}{%
\special{pn 4}%
\special{sh 1}%
\special{ar 5896 1692 16 16 0  6.28318530717959E+0000}%
\special{sh 1}%
\special{ar 5896 1692 16 16 0  6.28318530717959E+0000}%
}}%
%
{\color[named]{Black}{%
\special{pn 20}%
\special{pa 1696 2252}%
\special{pa 5896 1694}%
\special{fp}%
}}%
%
{\color[named]{Black}{%
\special{pn 8}%
\special{ar 3054 1718 98 354  1.5707963  4.7123890}%
}}%
%
{\color[named]{Black}{%
\special{pn 8}%
\special{pa 3054 1364}%
\special{pa 3059 1364}%
\special{pa 3059 1365}%
\special{pa 3061 1365}%
\special{fp}%
\special{pa 3086 1383}%
\special{pa 3086 1384}%
\special{pa 3087 1384}%
\special{pa 3087 1385}%
\special{pa 3089 1387}%
\special{pa 3089 1388}%
\special{pa 3090 1388}%
\special{pa 3090 1388}%
\special{fp}%
\special{pa 3105 1416}%
\special{pa 3106 1417}%
\special{pa 3106 1419}%
\special{pa 3107 1419}%
\special{pa 3107 1421}%
\special{pa 3108 1422}%
\special{pa 3108 1423}%
\special{fp}%
\special{pa 3119 1454}%
\special{pa 3119 1454}%
\special{pa 3120 1455}%
\special{pa 3120 1457}%
\special{pa 3121 1458}%
\special{pa 3121 1461}%
\special{fp}%
\special{pa 3130 1494}%
\special{pa 3130 1496}%
\special{pa 3131 1497}%
\special{pa 3131 1501}%
\special{pa 3131 1501}%
\special{fp}%
\special{pa 3138 1536}%
\special{pa 3138 1538}%
\special{pa 3139 1539}%
\special{pa 3139 1543}%
\special{fp}%
\special{pa 3144 1579}%
\special{pa 3144 1582}%
\special{pa 3145 1583}%
\special{pa 3145 1587}%
\special{fp}%
\special{pa 3149 1624}%
\special{pa 3149 1624}%
\special{pa 3149 1631}%
\special{fp}%
\special{pa 3151 1668}%
\special{pa 3151 1676}%
\special{fp}%
\special{pa 3152 1714}%
\special{pa 3152 1722}%
\special{fp}%
\special{pa 3151 1759}%
\special{pa 3151 1767}%
\special{fp}%
\special{pa 3149 1805}%
\special{pa 3149 1811}%
\special{pa 3148 1812}%
\special{pa 3148 1812}%
\special{fp}%
\special{pa 3145 1849}%
\special{pa 3145 1853}%
\special{pa 3144 1854}%
\special{pa 3144 1857}%
\special{fp}%
\special{pa 3139 1893}%
\special{pa 3139 1896}%
\special{pa 3138 1897}%
\special{pa 3138 1900}%
\special{fp}%
\special{pa 3131 1935}%
\special{pa 3131 1938}%
\special{pa 3130 1939}%
\special{pa 3130 1943}%
\special{fp}%
\special{pa 3121 1977}%
\special{pa 3120 1978}%
\special{pa 3120 1981}%
\special{pa 3119 1981}%
\special{pa 3119 1983}%
\special{fp}%
\special{pa 3107 2015}%
\special{pa 3107 2015}%
\special{pa 3107 2016}%
\special{pa 3106 2017}%
\special{pa 3106 2019}%
\special{pa 3105 2020}%
\special{pa 3105 2021}%
\special{pa 3104 2022}%
\special{fp}%
\special{pa 3088 2050}%
\special{pa 3087 2051}%
\special{pa 3087 2052}%
\special{pa 3086 2052}%
\special{pa 3086 2053}%
\special{pa 3085 2053}%
\special{pa 3085 2054}%
\special{pa 3084 2054}%
\special{pa 3084 2054}%
\special{fp}%
\special{pa 3061 2071}%
\special{pa 3059 2071}%
\special{pa 3059 2072}%
\special{pa 3054 2072}%
\special{fp}%
}}%
%
{\color[named]{Black}{%
\special{pn 8}%
\special{pa 4254 1500}%
\special{pa 4258 1500}%
\special{pa 4258 1501}%
\special{pa 4261 1501}%
\special{fp}%
\special{pa 4280 1521}%
\special{pa 4280 1523}%
\special{pa 4281 1523}%
\special{pa 4281 1525}%
\special{pa 4282 1525}%
\special{pa 4282 1527}%
\special{fp}%
\special{pa 4294 1556}%
\special{pa 4294 1559}%
\special{pa 4295 1559}%
\special{pa 4295 1563}%
\special{fp}%
\special{pa 4303 1595}%
\special{pa 4303 1599}%
\special{pa 4304 1600}%
\special{pa 4304 1602}%
\special{fp}%
\special{pa 4309 1637}%
\special{pa 4309 1645}%
\special{fp}%
\special{pa 4312 1681}%
\special{pa 4312 1681}%
\special{pa 4312 1689}%
\special{fp}%
\special{pa 4312 1727}%
\special{pa 4312 1735}%
\special{fp}%
\special{pa 4309 1771}%
\special{pa 4309 1779}%
\special{fp}%
\special{pa 4304 1813}%
\special{pa 4304 1816}%
\special{pa 4303 1817}%
\special{pa 4303 1820}%
\special{fp}%
\special{pa 4296 1853}%
\special{pa 4295 1854}%
\special{pa 4295 1857}%
\special{pa 4294 1857}%
\special{pa 4294 1860}%
\special{fp}%
\special{pa 4283 1889}%
\special{pa 4283 1889}%
\special{pa 4282 1889}%
\special{pa 4282 1891}%
\special{pa 4281 1891}%
\special{pa 4281 1893}%
\special{pa 4280 1893}%
\special{pa 4280 1894}%
\special{fp}%
\special{pa 4261 1915}%
\special{pa 4258 1915}%
\special{pa 4258 1916}%
\special{pa 4254 1916}%
\special{fp}%
}}%
%
{\color[named]{Black}{%
\special{pn 8}%
\special{ar 4254 1708 58 208  1.5707963  4.7123890}%
}}%
%
{\color[named]{Black}{%
\special{pn 20}%
\special{ar 5896 4330 94 362  1.5707963  4.7123890}%
}}%
%
{\color[named]{Black}{%
\special{pn 20}%
\special{pa 5896 3968}%
\special{pa 5901 3968}%
\special{pa 5901 3969}%
\special{pa 5903 3969}%
\special{fp}%
\special{pa 5926 3986}%
\special{pa 5926 3987}%
\special{pa 5927 3988}%
\special{pa 5927 3989}%
\special{pa 5928 3989}%
\special{pa 5928 3990}%
\special{pa 5929 3991}%
\special{pa 5929 3992}%
\special{pa 5929 3992}%
\special{fp}%
\special{pa 5943 4018}%
\special{pa 5944 4018}%
\special{pa 5944 4020}%
\special{pa 5945 4020}%
\special{pa 5945 4022}%
\special{pa 5946 4022}%
\special{pa 5946 4023}%
\special{fp}%
\special{pa 5956 4051}%
\special{pa 5956 4053}%
\special{pa 5957 4053}%
\special{pa 5957 4056}%
\special{pa 5958 4056}%
\special{pa 5958 4057}%
\special{fp}%
\special{pa 5966 4089}%
\special{pa 5966 4091}%
\special{pa 5967 4091}%
\special{pa 5967 4095}%
\special{pa 5967 4095}%
\special{fp}%
\special{pa 5974 4129}%
\special{pa 5974 4131}%
\special{pa 5975 4131}%
\special{pa 5975 4136}%
\special{fp}%
\special{pa 5980 4170}%
\special{pa 5980 4171}%
\special{pa 5981 4172}%
\special{pa 5981 4178}%
\special{fp}%
\special{pa 5985 4214}%
\special{pa 5985 4219}%
\special{pa 5986 4219}%
\special{pa 5986 4221}%
\special{fp}%
\special{pa 5988 4257}%
\special{pa 5988 4264}%
\special{fp}%
\special{pa 5990 4301}%
\special{pa 5990 4309}%
\special{fp}%
\special{pa 5990 4347}%
\special{pa 5990 4355}%
\special{fp}%
\special{pa 5989 4392}%
\special{pa 5989 4394}%
\special{pa 5988 4395}%
\special{pa 5988 4399}%
\special{fp}%
\special{pa 5986 4435}%
\special{pa 5986 4440}%
\special{pa 5985 4441}%
\special{pa 5985 4443}%
\special{fp}%
\special{pa 5982 4478}%
\special{pa 5982 4480}%
\special{pa 5981 4481}%
\special{pa 5981 4485}%
\special{fp}%
\special{pa 5976 4521}%
\special{pa 5976 4523}%
\special{pa 5975 4523}%
\special{pa 5975 4527}%
\special{fp}%
\special{pa 5968 4561}%
\special{pa 5968 4565}%
\special{pa 5967 4565}%
\special{pa 5967 4568}%
\special{fp}%
\special{pa 5958 4601}%
\special{pa 5958 4601}%
\special{pa 5958 4603}%
\special{pa 5957 4604}%
\special{pa 5957 4607}%
\special{pa 5956 4607}%
\special{pa 5956 4607}%
\special{fp}%
\special{pa 5946 4638}%
\special{pa 5946 4638}%
\special{pa 5945 4638}%
\special{pa 5945 4640}%
\special{pa 5944 4640}%
\special{pa 5944 4642}%
\special{pa 5943 4642}%
\special{pa 5943 4643}%
\special{fp}%
\special{pa 5930 4668}%
\special{pa 5929 4668}%
\special{pa 5929 4670}%
\special{pa 5928 4670}%
\special{pa 5928 4671}%
\special{pa 5927 4671}%
\special{pa 5927 4672}%
\special{pa 5926 4673}%
\special{fp}%
\special{pa 5903 4691}%
\special{pa 5901 4691}%
\special{pa 5901 4692}%
\special{pa 5896 4692}%
\special{fp}%
}}%
%
{\color[named]{Black}{%
\special{pn 20}%
\special{pa 1694 4200}%
\special{pa 5894 3966}%
\special{fp}%
}}%
%
{\color[named]{Black}{%
\special{pn 4}%
\special{sh 1}%
\special{ar 5896 4692 16 16 0  6.28318530717959E+0000}%
\special{sh 1}%
\special{ar 5896 4692 16 16 0  6.28318530717959E+0000}%
}}%
%
{\color[named]{Black}{%
\special{pn 20}%
\special{pa 5896 4692}%
\special{pa 1696 4204}%
\special{fp}%
}}%
%
{\color[named]{Black}{%
\special{pn 4}%
\special{sh 1}%
\special{ar 1696 4204 8 8 0  6.28318530717959E+0000}%
\special{sh 1}%
\special{ar 1696 4204 8 8 0  6.28318530717959E+0000}%
}}%
%
{\color[named]{Black}{%
\special{pn 4}%
\special{sh 1}%
\special{ar 1696 4204 16 16 0  6.28318530717959E+0000}%
\special{sh 1}%
\special{ar 1696 4204 16 16 0  6.28318530717959E+0000}%
}}%
%
{\color[named]{Black}{%
\special{pn 8}%
\special{ar 3056 4244 34 120  1.5707963  4.7123890}%
}}%
%
{\color[named]{Black}{%
\special{pn 8}%
\special{pa 3056 4124}%
\special{pa 3059 4124}%
\special{pa 3059 4125}%
\special{pa 3061 4125}%
\special{pa 3061 4126}%
\special{pa 3062 4126}%
\special{fp}%
\special{pa 3078 4151}%
\special{pa 3078 4154}%
\special{pa 3079 4154}%
\special{pa 3079 4157}%
\special{pa 3080 4157}%
\special{pa 3080 4158}%
\special{fp}%
\special{pa 3087 4193}%
\special{pa 3087 4199}%
\special{pa 3088 4199}%
\special{pa 3088 4200}%
\special{fp}%
\special{pa 3090 4240}%
\special{pa 3090 4248}%
\special{fp}%
\special{pa 3088 4287}%
\special{pa 3088 4289}%
\special{pa 3087 4289}%
\special{pa 3087 4295}%
\special{fp}%
\special{pa 3080 4330}%
\special{pa 3080 4331}%
\special{pa 3079 4331}%
\special{pa 3079 4334}%
\special{pa 3078 4334}%
\special{pa 3078 4337}%
\special{fp}%
\special{pa 3062 4362}%
\special{pa 3061 4362}%
\special{pa 3061 4363}%
\special{pa 3059 4363}%
\special{pa 3059 4364}%
\special{pa 3056 4364}%
\special{fp}%
}}%
%
{\color[named]{Black}{%
\special{pn 8}%
\special{pa 4256 4058}%
\special{pa 4259 4058}%
\special{pa 4259 4059}%
\special{pa 4262 4059}%
\special{pa 4262 4060}%
\special{pa 4262 4060}%
\special{fp}%
\special{pa 4279 4084}%
\special{pa 4280 4084}%
\special{pa 4280 4086}%
\special{pa 4281 4087}%
\special{pa 4281 4089}%
\special{pa 4282 4089}%
\special{pa 4282 4090}%
\special{fp}%
\special{pa 4291 4121}%
\special{pa 4291 4123}%
\special{pa 4292 4124}%
\special{pa 4292 4128}%
\special{pa 4293 4129}%
\special{fp}%
\special{pa 4299 4163}%
\special{pa 4299 4163}%
\special{pa 4299 4170}%
\special{pa 4300 4171}%
\special{fp}%
\special{pa 4303 4208}%
\special{pa 4303 4210}%
\special{pa 4304 4211}%
\special{pa 4304 4215}%
\special{fp}%
\special{pa 4306 4253}%
\special{pa 4306 4261}%
\special{fp}%
\special{pa 4306 4300}%
\special{pa 4306 4308}%
\special{fp}%
\special{pa 4304 4345}%
\special{pa 4304 4349}%
\special{pa 4303 4350}%
\special{pa 4303 4353}%
\special{fp}%
\special{pa 4299 4389}%
\special{pa 4299 4389}%
\special{pa 4299 4397}%
\special{fp}%
\special{pa 4293 4430}%
\special{pa 4293 4432}%
\special{pa 4292 4432}%
\special{pa 4292 4436}%
\special{pa 4291 4436}%
\special{pa 4291 4436}%
\special{fp}%
\special{pa 4282 4469}%
\special{pa 4282 4471}%
\special{pa 4281 4471}%
\special{pa 4281 4473}%
\special{pa 4280 4474}%
\special{pa 4280 4475}%
\special{fp}%
\special{pa 4262 4500}%
\special{pa 4262 4500}%
\special{pa 4262 4501}%
\special{pa 4259 4501}%
\special{pa 4259 4502}%
\special{pa 4256 4502}%
\special{fp}%
}}%
%
{\color[named]{Black}{%
\special{pn 8}%
\special{ar 4256 4280 50 222  1.5707963  4.7123890}%
}}%
%
{\color[named]{Black}{%
\special{pn 8}%
\special{pa 5490 3200}%
\special{pa 2102 3200}%
\special{fp}%
}}%
\put(16.9000,-32.0000){\makebox(0,0){$t=0$}}%
\put(58.9000,-32.0000){\makebox(0,0){$t=1$}}%
\put(37.0000,-9.0000){\makebox(0,0){Figure of $\Delta $}}%
\put(37.0000,-38.0000){\makebox(0,0){Figure of $\Gamma  $}}%
\end{picture}%

\end{center}


\vspace{0.5cm}

Let $r_1(t)<r_2(t)<r_3(t)$ be the real roots
of $x^3+(3x+4t^l)^2$ for $0<t<1$. Then one has
\begin{equation}\label{Example-4}
\int_{\Delta} t^{j-1}dt\frac{dx}{y}=
2\sqrt{-3}
\overbrace{\int_0^1 t^{j-1}dt\int_{r_1(t)}^{r_2(t)}
\frac{dx}{\sqrt{x^3+(3x+4t^l)^2}}}^{\mbox{positive real number}}
\in \sqrt{-1}\R_{>0}
\end{equation}
\begin{equation}\label{Example-5}
\int_\Gamma t^{j-1}dt\frac{dx}{y}=
2\sqrt{-3}
\overbrace{\int_0^1 t^{j-1}dt\int_{r_2(t)}^{r_3(t)}
\frac{dx}{\sqrt{x^3+(3x+4t^l)^2}}}^{\mbox{purely imaginary number}\ne0}
\in \R_{>0}.
\end{equation}
Let $\sigma:X_\C\to X_\C$ be an automorphism given by $\sigma(x,y,t)=(x,y,\zeta_l t)$.
Since 
\[\sigma^*t^{j-1}dt\frac{dx}{y}=\zeta^j_l t^{j-1}dt\frac{dx}{y},\] 
one has
\[
\int_{\sigma_*^k\Delta} t^{j-1}dt\frac{dx}{y}=
\zeta_l^{kj}\int_{\Delta} t^{j-1}dt\frac{dx}{y}.
\]
This and an elementary calculation show
that 
\[
\Delta-\sigma_*\Delta,\Delta-\sigma^2_*\Delta,
\cdots,\Delta-\sigma^{l-1}_*\Delta
\]
are linearly independent in $\bE(\ol{U},\Q)$.
Since it is $(l-1)$-dimensional by \eqref{Example-7}, the above is a basis of $\bE(\ol{U},\Q)$.
$\bE(\ol{U},D;\Q)/\bE(\ol{U},\Q)$ is $l$-dimensional with a basis
\[
\Gamma,\sigma_*\Gamma,
\cdots,\sigma^{l-1}_*\Gamma.
\]
Let $F_\infty$ denotes the infinite Frobenius morphism.
Then 
\[
F_\infty(\Delta)=-\Delta,\quad F_\infty(\Gamma)=\Gamma.
\]
By Thm.\ref{ExtMHS} and the above computations
we have the following.
\begin{thm}\label{Example-thm0}
Suppose $(l,6)=1$. Put $h=\dim F^1 V_\dR=l-\lfloor\frac{l-1}{3}\rfloor-1$ and
$\zeta=\exp(2\pi i/l)$.
Let
\[
A=
\left((\zeta^{pq}-\zeta^{-pq})
\int_{\Delta} t^{p-1}dt\frac{dx}{y}\right)_{1\leq p\leq h,~1\leq q \leq (l-1)/2}
\]
be $h\times (l-1)/2$-matrix (the entries are real numbers by \eqref{Example-4}).
Then
\[
\mathrm{Ext}^1_{\R\mbox{-}\MHS}(\R,H^2(X)_\ind\ot\R(1))^{F_\infty=1}
\cong\Coker[A:\R^{(l-1)/2}\lra \R^h].
\]
and we have
\[
\ol{\reg}_\R(\xi_{D_1})=
\pm\left(
\int_{\Gamma} dt\frac{dx}{y},
\cdots,
\int_{\Gamma} t^{h-1}dt\frac{dx}{y}
\right)\in \R^h/\Image A
\]
under the above isomorphism.
\end{thm}
\begin{cor}\label{Example-cor}
Suppose $(l,6)=1$.
Then
\[
\ol{\reg}_\R(\xi_{D_1})\ne0\in\Ext^1_{\MHS}(\R,H^2(X)_\ind\ot\R(1))^{F_\infty=1}.
\]
In particular $\xi_{D_1}$ is regulator indecomposable.
\end{cor}
\begin{pf}
Put $h:=l-\lfloor\frac{l-1}{3}\rfloor-1$ and
$\zeta:=\exp(2\pi i/l)$.
Put
\[
I_p:=\int_\Delta t^{p-1}dt\frac{dx}{y},
\quad
J_p:=\int_\Gamma t^{p-1}dt\frac{dx}{y}.
\]
Then
\[
\ol{\reg}_\R(\xi_{D_1})\ne0\in\Ext^1_{\MHS}(\R,H^2(X)_\ind\ot\R(1))^{F_\infty=1}
\]
if and only if the rank of a matrix
\begin{equation}\label{Example-cor-1}
\begin{pmatrix}
 (\zeta-\zeta^{-1})I_{1}& (\zeta^2-\zeta^{-2})I_{1}
&\cdots&(\zeta^{(l-1)/2}-\zeta^{-(l-1)/2})I_{1}&J_1\\
 (\zeta^2-\zeta^{-2})I_{2}& (\zeta^{4}-\zeta^{-4})I_{2}
&\cdots&(\zeta^{l-1}-\zeta^{-(l-1)})I_{2}&J_2\\
 \vdots&\vdots&&\vdots&\vdots\\
(\zeta^{h}-\zeta^{-h})I_{h}& (\zeta^{2h}-\zeta^{-2h})I_{h}
&\cdots&(\zeta^{h(l-1)/2}-\zeta^{-h(l-1)/2})I_{h}&J_h\\
\end{pmatrix}
\end{equation}
is maximal.
Thus it is enough to show that
\begin{equation}\label{Example-cor-2}
\det
\begin{pmatrix}
 (\zeta-\zeta^{-1})& (\zeta^2-\zeta^{-2})
&\cdots&(\zeta^{(l-1)/2}-\zeta^{-(l-1)/2})&J_1/I_1\\
 (\zeta^2-\zeta^{-2})& (\zeta^{4}-\zeta^{-4})
&\cdots&(\zeta^{l-1}-\zeta^{-(l-1)})&J_2/I_2\\
 \vdots&\vdots&&\vdots&\vdots\\
(\zeta^{k}-\zeta^{-k})& (\zeta^{2k}-\zeta^{-2k})
&\cdots&(\zeta^{k(l-1)/2}-\zeta^{-k(l-1)/2})&J_{k}/I_{k}\\
\end{pmatrix}
\end{equation}
is nonzero where $k=(l+1)/2$.
Since the sum of the $(k-1)$-th row and $k$-th row is
$
(0,\cdots,0,J_{k-1}/I_{k-1}+J_{k}/I_k)
$,
one has
\begin{align*}
\eqref{Example-cor-2}=&
(J_{k-1}/I_{k-1}+J_{k}/I_k)\times\det(\zeta^{pq}-\zeta^{-pq})_{1\leq p,q\leq (l-1)/2}\\
=&(J_{k-1}/I_{k-1}+J_{k}/I_k)\times\sqrt{(-l)^{(l-1)/2}}.
\end{align*}
Since $J_p/I_p\in i\R_{>0}$ by \eqref{Example-4} and \eqref{Example-5}, 
this is non-zero. 
\end{pf}

\subsection{Another description of $\int_{\Delta}$ and $\int_{\Gamma}$}
When $l=1$, $f:X\to \P^1$ is the universal elliptic curve over $X_1(3)$.
Using this, one can obtain another description of the real regulator.

\medskip

Let $q=\exp(2\pi i z)$ and
\[
E_{3a}(z):=1-9\sum_{n=1}^\infty
\left(\sum_{k|n}\left(\frac{k}{3}\right)k^2\right)q^n,
\]
\[
E_{3b}(z):=\sum_{n=1}^\infty\left(\sum_{k|n}\left(\frac{n/k}{3}\right)k^2\right)q^n
\]
be the Eisenstein series of weight 3 for $\Gamma_1(3)$, where $(\frac{k}{3})$ denotes 
the Legendre symbol.
Then
\[
t^l
=\frac{E_{3a}}{E_{3a}+27E_{3b}}
\]
and
\[
l\frac{dt}{t}\frac{dx}{y}=-27E_{3b}\frac{du}{u}\frac{dq}{q},\quad
\frac{lt^{l-1}dt}{t^l-1}\frac{dx}{y}=E_{3a}\frac{du}{u}\frac{dq}{q}
\]
where ``$du/u$" denotes the canonical invariant 1-form of the Tate curve around 
the cusp $z=i\infty$ ($t=1$).
Therefore we have
\begin{equation}\label{Example-1}
\int_{\Delta} t^{j-1}dt\frac{dx}{y}=\frac{-27}{l}\times (2\pi i)^2
\int_0^{i\infty} t^{j}E_{3b}(z)dz
\end{equation}
\begin{equation}\label{Example-2}
\int_\Gamma t^{j-1}dt\frac{dx}{y}=\frac{-27}{l}\times (2\pi i)^2
\int_0^{i\infty} t^{j}E_{3b}(z)zdz.
\end{equation}
On the other hand there are formulas
\begin{equation}\label{Example-3}
\frac{E_{3a}}{E_{3a}+27E_{3b}}(\frac{-1}{3z})
=\frac{27E_{3b}}{E_{3a}+27E_{3b}}(z),\quad
27E_{3b}(\frac{-1}{3z})=
3\sqrt{3}iz^3 E_{3a}(z)
\end{equation}
on the Eisenstein series.
Applying \eqref{Example-3} to \eqref{Example-1} and \eqref{Example-2},
we have the following theorem.

\begin{thm}\label{Example-thm}
Put $c:=\exp(-2\pi/\sqrt{3})=0.026579933\cdots$. 
Define rational numbers $a_n(j)$ and $b_n(j)$ by
\begin{multline*}
E_{3b}\left(\frac{E_{3a}}{E_{3a}+27E_{3b}}\right)^{j/l}=\sum_{n=1}^\infty a_n{(j)}q^n\\
=q+\left(3-27\frac{j}{l}\right)q^2
+\left(9-\frac{81}{2}\frac{j}{l}+\frac{729}{2}\left(\frac{j}{l}\right)^2\right)q^3+\cdots,
\end{multline*}
\begin{multline*}
E_{3a}\left(\frac{E_{3b}}{q(E_{3a}+27E_{3b})}\right)^{j/l}
=\sum_{n=0}^\infty b_n{(j)}q^n\\
=1+\left(-9-15\frac{j}{l}\right)q+
\left(27+\frac{387}{2}\frac{j}{l}+\frac{225}{2}\left(\frac{j}{l}\right)^2\right)q^2+\cdots.
\end{multline*}
Put
\begin{align*}
I(j)&=\sum_{n=1}^\infty \frac{a_n{(j)}}{n}c^n+3^{3j/l-3}\sum_{n=0}^\infty b_n{(j)}\left(\frac{1}{n+j/l}
+\frac{\sqrt{3}}{2\pi(n+j/l)^2}\right)c^{n+j/l}\\
J(j)&=\sum_{n=1}^\infty a_n(j)\left(
\frac{2\pi}{\sqrt{3}n}+\frac{1}{n^2}\right)c^n+2\pi\cdot 3^{3j/l-7/2}\sum_{n=0}^\infty 
\frac{b_n(j)}{n+j/l}c^{n+j/l}.
\end{align*}
Then we have
\[
\int_{\Delta} t^{j-1}dt\frac{dx}{y}
=\frac{54\pi i}{l}I(j),\quad 
\int_\Gamma t^{j-1}dt\frac{dx}{y}
=\frac{-27}{l}J(j)
\]
for $1\leq j\leq l-1$.
\end{thm}
This is useful since the series $I(j)$ and $J(j)$ converge rapidly !

\begin{exmp}\label{exp-l5}{\rm Suppose $l=5$. Then $X$ is a K3 surface.
By Thm.\ref{Example-thm}, one has 
\begin{center}
\begin{tabular}{|c|c|c|}
\hline
&$I(j)$
&$J(j)$\\
\hline
$j=1$&
$ 0.42745977255318$&$0.717696894965804$
\\
\hline
$j=2$&
$0.151180954233147$
&$0.377159120670032$
\\
\hline
$j=3$&
$0.0871841692346256$
&$0.261572572611421$
\\
\hline
$j=4$&
$0.0603840144077692$&$0.202670503662525$
\\
\hline
\end{tabular}
\end{center}
\[
\mathrm{Ext}^1_{\R\mbox{-}\MHS}(\R,H^2(X)_\ind\ot\R(1))^{F_\infty=1}
\cong\Coker(\R^2\os{A}{\lra} \R^3).
\]
Since this is 1-dimensional,
this has the canonical base $e_{\ind,\Q}$ (up to $\Q^\times$) and a {\it different} base
$e_{\ind,\Q}^\ff$
(\S \ref{false-sect}).
With respect to $e_{\ind,\Q}^\ff$, one has
\[
\ol{\reg}_\R(\xi_{D_1})=\pi^2
\begin{vmatrix}
 i(\zeta-\zeta^{-1})I(1)& i(\zeta^2-\zeta^{-2})I(1)&J(1)\\
 i(\zeta^2-\zeta^{-2})I(2)& i(\zeta^4-\zeta^{-4})I(2)&J(2)\\
i(\zeta^3-\zeta^{-3})I(3)& i(\zeta^6-\zeta^{-6})I(3)&J(3)\\
\end{vmatrix}\mod \Q^\times \quad(\zeta:=\exp(2\pi i/5)).
\]
Since $s=(l-1)/2=2$ and $\det H_\dR^2(X/\Q)_\ind\ot
[\det H_B^2(X_\C)_\ind]^{-1}=\sqrt{5}$ 
one has 
\begin{align*}
\ol{\reg}_\R(\xi_{D_1})
&=\frac{\sqrt{5}}{\pi^2}\cdot\pi^2
\begin{vmatrix}
 i(\zeta-\zeta^{-1})I(1)& i(\zeta^2-\zeta^{-2})I(1)&J(1)\\
 i(\zeta^2-\zeta^{-2})I(2)& i(\zeta^4-\zeta^{-4})I(2)&J(2)\\
i(\zeta^3-\zeta^{-3})I(3)& i(\zeta^6-\zeta^{-6})I(3)&J(3)\\
\end{vmatrix}\\
&=-5\sqrt{5}I(1)I(2)I(3)\left(\frac{J(2)}{I(2)}+\frac{J(3)}{I(3)}\right)\\
&=0.346139631939354
\mod \Q^\times 
\end{align*}
with respect to $e_{\ind,\Q}$ by Prop.\ref{qstr-8}.
}
\end{exmp}

\begin{exmp}\label{exp-l7} {\rm Suppose $l=7$. Then $h^{20}(X)=h^{02}(X)=2$, $h^{11}(X)=30$.
\begin{center}
\begin{tabular}{|c|c|c|}
\hline
&$I(j)$
&$J(j)$\\
\hline
$j=1$&
$0.740059830730164$&$0.987994510350351$
\\
\hline
$j=2$&
$0.24646699651114$
&$0.51401702238944$
\\
\hline
$j=3$&
$0.137265313181901$
&$0.354195498081428$
\\
\hline
$j=4$&
$0.0929578147374374$&$0.273237679671921$
\\
\hline
$j=5$&
$0.0696363855176379$&$0.224004116344261$
\\
\hline
$j=6$&
$0.0554349861351089$&$0.19073921727221$
\\
\hline
\end{tabular}
\end{center}
\[
\mathrm{Ext}^1_{\R\mbox{-}\mathrm{MHS}}(\R,H^2(X)_\ind\ot\R(1))^{F_\infty=1}
\cong\Coker(\R^3\os{A}{\lra} \R^4).
\]
Since $s=(l-1)/2=3$ and $\det H_\dR^2(X/\Q)_\ind\ot
[\det H_B^2(X_\C)_\ind]^{-1}=\sqrt{-7}$,
one has 
\begin{align*}
\ol{\reg}_\R(\xi_{D_1})&=\frac{\sqrt{7}}{\pi^3}\cdot\pi^3
\begin{vmatrix}
 i(\zeta-\zeta^{-1})I(1)& i(\zeta^2-\zeta^{-2})I(1)& i(\zeta^3-\zeta^{-3})I(1)&J(1)\\
 i(\zeta^2-\zeta^{-2})I(2)& i(\zeta^4-\zeta^{-4})I(2)& i(\zeta^6-\zeta^{-6})I(2)&J(2)\\
i(\zeta^3-\zeta^{-3})I(3)& i(\zeta^6-\zeta^{-6})I(3)& i(\zeta^9-\zeta^{-9})I(3)&J(3)\\
i(\zeta^4-\zeta^{-4})I(4)& i(\zeta^8-\zeta^{-8})I(4)& i(\zeta^{12}-\zeta^{-12})I(4)&J(4)\\
\end{vmatrix}\\
&=49I(1)I(2)I(3)I(4)\left(\frac{J(3)}{I(3)}+\frac{J(4)}{I(4)}\right)\\
&=0.629487860860585
\mod \Q^\times\quad(\zeta:=\exp(2\pi i/7))
\end{align*}
with respect to the canonical $\Q$-structure $e_{\ind,\Q}$.
}
\end{exmp}
\begin{rem}
According to the Beilinson conjecture, $\ol{\reg}_\R(\xi_{D_1})$
in Example \ref{exp-l5} or \ref{exp-l7} is expected to be the value of 
the $L$-function $L(h^2(X)_\ind,s)$ at $s=1$ $(\cite{schneider})$. 
\end{rem}

\def\EE{2g_2g'_3-3g'_2g_3}
\def\EEE{6g_2g'_3-9g'_2g_3}

\section{Appendix : Gauss-Manin connection}\label{Appendix-sect}
Let $R$ be an integral domain of characteristic $0$ in which $6$ is invertible.
For a smooth scheme $Y$ over $T$, we denote by $\Omega^q_{Y/T}=\os{q}{\wedge}_{\O_Y}
\Omega^1_{Y/T}$ the sheaf of relative differential $q$-forms on $Y$ over $T$. 
If $T=\Spec R$, we
simply write $\Omega^q_{Y}=\Omega^q_{Y/R}$.

\subsection{Explicit formulas}\label{expPF-sect}
Let $S$ be an irreducible affine smooth scheme over $R$ of relative dimension one.
Let $g_2,g_3\in \O_S(S)=\vg(S,\O_S)$ satisfy $\Delta:=g_2^3-27g_3^2\in \O_S(S)^\times$.
Let $f:U\to S$ be a projective smooth family of
elliptic curves whose affine form is given by a Weierstrass equation
$y^2=4x^3-g_2x-g_3$.
More precisely letting 
\[
U_0=\Spec \O_S(S)[x,y]/(y^2-4x^3+g_2x+g_3),
\]
\[
U_\infty=\Spec \O_S(S)[u,z]/(z-4u^3+g_2uz^2+g_3z^3),
\]
$U$ is obtained by gluing
$U_0$ and $U_\infty$ via identification 
$u=x/y,~z=1/y$.
Let $e:S\lra U$ be a section given by $(u,z)=(0,0)$.
To describe the de Rham cohomology
$H^q_\dR(U/S):=\check{H}^q(U,\Omega^\bullet_{U/S})$
we use the Cech complex.
Write
\[
\check{C}^0({\mathscr F}):=\vg(U_0,{\mathscr F})\op\vg(U_\infty,{\mathscr F}),\quad
\check{C}^1({\mathscr F}):=\vg(U_0\cap U_\infty,{\mathscr F})
\]
for a (Zariski) sheaf $\mathscr F$.
Then the double complex
\[
\xymatrix{
\check{C}^0(\O_U)\ar[r]^d\ar[d]_\delta&
\check{C}^0(\Omega^1_{U/S})\ar[d]^\delta&(x_0,x_\infty)\ar[d]^\delta
\ar[r]^d&(dx_0,dx_\infty)\\
\check{C}^1(\O_U)\ar[r]^d&
\check{C}^1(\Omega^1_{U/S})&x_0-x_\infty
}
\]
gives rise to the total complex
\[
\check{C}^\bullet(U/S):\check{C}^0(\O_U)\os{\delta\times d}{\lra} 
\check{C}^1(\O_U)\times \check{C}^0(\Omega^1_{U/S})
\os{(-d)\times\delta}{\lra} \check{C}^1(\Omega^1_{U/S})
\]
of $R$-modules starting from degree 0, and the cohomology of it is the de Rham cohomology
$H^\bullet_\dR(U/S)$:
\[
H^q_\dR(U/S)=H^q(\check{C}^\bullet(U/S)),\quad q\geq 0.
\]
Elements of $H^1_\dR(U/S)$ are represented by cocycles
\[
(f)\times (x_0,x_\infty)\quad \mbox{with }df=x_0-x_\infty.
\]

The purpose of Appendix is to write down
the {\it Gauss-Manin connection}
\[
\nabla:H^1_\dR(U/S)\lra \Omega^1_S\ot
H^1_\dR(U/S)
\]
(we use the same symbol ``$\Omega^1_S$" for $\vg(S,\Omega^1_S)$ 
since it will be clear from the context which is meant).
This is defined in the following way (cf. \cite{h} Ch.III, \S 4).
By applying $Rf_*$ on an exact sequence
\[
0\lra
f^*\Omega^1_S\ot\Omega^{\bullet-1}_{U/S}
\lra
\Omega^\bullet_U\lra \Omega^\bullet_{U/S}\lra 0,
\]
one has the connecting homomorphism
$R^1f_*\Omega^\bullet_{U/S}\to R^2f_*(f^*\Omega^1_S\ot
\Omega^{\bullet-1}_{U/S})\cong \Omega^1_S\ot R^2f_*(
\Omega^{\bullet-1}_{U/S})$.
By identifying $R^2f_*(\Omega^{\bullet-1}_{U/S})$ with
$R^1f_*\Omega^\bullet_{U/S}$, one gets the Gauss-Manin connection $\nabla$.
Here we should be careful about ``sign" because
the differential of the complex $\Omega^{\bullet-1}_{U/S}$ is ``$-d$" :
\[
\Omega^{\bullet-1}_{U/S}:\O_U\os{-d}{\lra} \Omega^1_{U/S}
\]
where the first term is placed in degree 1.
So we need to choose an isomorphism between $R^qf_*\Omega^{\bullet}_{U/S}$
and $R^{q+1}f_*\Omega^{\bullet-1}_{U/S}$
because the natural one is unique up to sign.
Here we choose it by 
\[
\xymatrix{
\O_U\ar[r]^d\ar[d]_{-\id}&\Omega^1_{U/S}\ar[d]^{\id}\\
\O_U\ar[r]^{-d}&\Omega^1_{U/S}
}
\]
Then $\nabla$ satisfies the usual Leibniz rule
\[
\nabla(fe)=df\ot e+f\nabla(e),\quad e\in H^1_\dR(U/S),~f\in\O(S).
\]
\begin{thm}\label{A1}
Suppose that $\Omega^1_S$ is a free $\O_S$-module with a base 
$dt\in \vg(S,\Omega^1_{S})$.
For $f\in \O_S(S)$, we define $f'\in\O_S(S)$ by $df=f'dt$.
Let
\begin{equation}\label{A1-1}
\canh:=(0)\times (\frac{dx}{y},\frac{dx}{y})
\end{equation}
\begin{equation}\label{A1-2}
\can:=(\frac{x^2}{y})\times (\frac{xdx}{2y},\frac{(2g_2x^2+3g_3x)dx}{2y^3})
\end{equation}
be elements in $H^1_\dR(U/S)$. 
Then we have
\begin{equation}\label{A1-3}
\nabla\left(\canh\right)=\left(\frac{\EEE}{\Delta}dt\ot\can
-\frac{\Delta^\prime}{12\Delta}dt\ot\canh\right)\in \Omega^1_S \ot H^1_\dR(U/S),
\end{equation}
\begin{equation}\label{A1-4}
\nabla\left(\can\right)=\left(\frac{\Delta^\prime}{12\Delta}dt\ot\can
-\frac{g_2(\EE)}{16\Delta}dt\ot\canh\right)\in 
\Omega^1_S \ot H^1_\dR(U/S).
\end{equation}
\end{thm}
Note that $\canh$ and $\can$ are basis of the free $\O(S)$-module
$H^1_\dR(U/S)$:
\[
H^1_\dR(U/S)=\O_S(S)\canh\op\O_S(S)\can.
\]
The following is straightforward from Thm. \ref{A1}.
\begin{cor}\label{PF-rem}
Put $\cH:=R^1f_*\Omega^\bullet_{U/S}$ and
\[
\cH^{1,0}:=f_*\Omega^1_{U/S}=\O_S\canh,\quad
\cH^{0,1}:=\cH/\cH^{1,0}\cong \O_S\can.
\]
Then the $\O_S$-linear map 
\begin{equation}\label{pf0}
\ol{\nabla}:\cH^{1,0}\lra \Omega^1_S\ot\cH^{0,1}
\end{equation}
induced from the Gauss-Manin connection $\nabla$ is described as follows.
\[
\ol{\nabla}\left(\canh\right)=\frac{\EEE}{\Delta}dt\ot\can.
\]
In particular, noting
\[
\frac{j'}{j}=27\cdot\frac{g_3}{g_2}\cdot\frac{\EE}{\Delta},\quad
j:=\frac{1728g_2^3}{g_2^3-27g_3^2},
\]
\eqref{pf0} is bijective if and only if
\[
\frac{g_2}{g_3}\frac{dj}{j}\in\Omega^1_S
\]
is a base of $\O_S$-module.
\end{cor}
Let us consider a diagram
\begin{equation}\label{pf1}
\xymatrix{
&0\ar[d]\\
&\Omega^1_S\ot\cH^{1,0}\ar[d]\\
\cH\ar[r]^{\nabla\quad}\ar[d]_=&\Omega^1_S\ot\cH\ar[d]\\
\cH\ar[r]^{\wt{\nabla}\quad}&\Omega^1_S\ot\cH^{0,1}\ar[d]\\
&0}
\end{equation}
Let $S^o\subset S$ be a Zariski open set such that 
$\ol{\nabla}$ is bijective on $S^o$.
Then it gives rise to an exact sequence
\begin{equation}\label{pf2}
\ker\wt{\nabla}|_{S^o}\lra \Omega^1_{S^o}\ot\cH^{1,0}|_{S^o}\lra \Coker\nabla|_{S^o} \lra0.
\end{equation}
Since the natural map $\ker\wt{\nabla}|_{S^o}\to \cH^{0,1}|_{S^o}$ is bijective,
we have an exact sequence
\begin{equation}\label{pf3}
\begin{CD}
0@>>>\vg(S^o,\cH^{0,1})@>{\PF}>>\vg(S^o,\Omega^1_{S}\ot\cH^{1,0})
@>>>
\vg(S^o,\Coker\nabla)@>>>0\\
&&&&&&@|\\
&&&&&&H^1_\dR(S^o,\cH).
\end{CD}
\end{equation}
The map $\PF$ in \eqref{pf3} is called the {\it Picard-Fuchs operator}.
\begin{cor}\label{PF}
Suppose that $\Omega^1_{S^o}$ is a free $\O_{S^o}$-module with a base 
$dt\in \vg(S^o,\Omega^1_{S^o})$.
Write $f'dt=df$ for $f\in \O(S^o)$.
Put 
\[
A:=\frac{-\Delta}{\EEE}
\]
\[
B:=\frac{1}{48}\left(
\frac{g_2(g'_2)^2-12(g'_3)^2}{\EE}
-\left(
\frac{4\Delta'}{3(\EE)}
\right)^\prime\right).
\]
Then the Picard-Fuchs operator is described as follows.
\[
\PF(f(t)\can)=(f^{\prime\prime}A+f'(t)A'+fB)dt \ot\canh,\quad f\in \O(S^o)
\]
\end{cor}
\begin{pf}
Let
\[
z:=f(t)\can-\frac{\Delta}{\EEE}(f'(t)+\frac{\Delta'}{12\Delta}f(t))
\canh\in \vg(S^o,\cH).
\]
This belongs to the kernel of $\wt{\nabla}$ by \eqref{A1-3} and \eqref{A1-4}.
Then 
\[
\PF(f(t)\can)=\nabla(z)
\]
and apply \eqref{A1-3} and \eqref{A1-4} again to the RHS.
\end{pf}

\subsection{Proof of Theorem \ref{A1}}
\begin{lem}\label{A0}
\[
(\frac{x^i}{y^j})\times(0,0)\equiv (0)\times 
(0,d(\frac{x^i}{y^j}))\quad (0\leq i\leq j)
\]
\[
(x^iy^j)\times(0,0)\equiv (0)\times 
(-d(x^iy^j),0)\quad (i,~j\geq 0)
\]
in $\check{C}^1(\O_U)\times
\check{C}^0(\Omega^1_{U/S})$
where ``$\equiv$" denote modulo $\Image\check{C}^0(\O_U)$.
\end{lem}
\begin{pf}
Straightforward from the definition.
\end{pf}
\begin{lem}\label{A2}
Let $\eta_U$ be the generic point of $U$. We think of $dx$ and $dt$ as elements in 
$\vg(\eta_U,\Omega^1_U)$.
Then 
\begin{equation}\label{A2-1}
\frac{dx}{y}\in \vg(U_\infty,\Omega^1_U),
\end{equation}
\begin{equation}\label{A2-2}
\widehat{\frac{dx}{y}}:=\frac{dx}{y}-\frac{(6g_2x^2-9g_3x-g_2^2)(g'_2x+g'_3)}{\Delta}\frac{dt}{y}
\in \vg(U_0,\Omega^1_U).
\end{equation}
\end{lem}
\begin{pf} 
\begin{align*}
\frac{dx}{y}&=du-\frac{u}{z}
((12u^2-g_2z^2)du-(2g_2uz+3g_3z^2)dz-(g'_2uz^2+g'_3z^3)dt)\\
&\equiv -\frac{12u^3}{z}du\mod \vg(U_\infty,\Omega^1_U)\\
&=-3(1+g_2uz+g_3z^2)du\\
&\equiv 0\mod \vg(U_\infty,\Omega^1_U).
\end{align*}
Hence \eqref{A2-1} follows.
Next we show \eqref{A2-2}.
Since
$f(x)=4x^3-g_2x-g_3$ is prime to $f'(x)=12x^2-g_2$, there are $a(x)$ and $b(x)$ such that
\[
a(x)f(x)+b(x)f'(x)=1.
\]
Explicitly, they are given as follows.
\[
a(x)=\frac{9(3g_3-2g_2x)}{\Delta},\quad
b(x)=\frac{6g_2x^2-9g_3x-g_2^2}{\Delta}.
\]
Now
\begin{align*}
\frac{dx}{y}&=\frac{a(x)f(x)dx+b(x)f'(x)dx}{y}\\
&=\frac{a(x)y^2dx+b(x)(2ydy+(g'_2x+g'_3)dt)}{y}\\
&=a(x)ydx+2b(x)dy+\frac{b(x)(g'_2x+g'_3)}{y}dt\in\vg(U_0\cap U_\infty,\Omega^1_U)
\end{align*}
and hence we have
\[
\widehat{\frac{dx}{y}}=\frac{dx}{y}-\frac{b(x)(g'_2x+g'_3)}{y}dt=a(x)ydx+2b(x)dy\in\vg(U_0,\Omega^1_U).
\]
\end{pf}

Let us prove \eqref{A1-3}.
Let
\[
\widehat{\canh}:=(0)\times (\widehat{\frac{dx}{y}},\frac{dx}{y})
\in \check{C}^1(\O_U)\times \check{C}^0(\Omega^1_{U})
\]
be a lifting of $\canh$ where $\widehat{\frac{dx}{y}}$ is as in \eqref{A2-2}.
Let 
\[
\xymatrix{
\check{C}^0(\O_U)\ar[r]^d\ar[d]_\delta&
\check{C}^0(\Omega^1_{U})\ar[d]^\delta\ar[r]^d&
\check{C}^0(\Omega^2_{U})\ar[d]^\delta\\
\check{C}^1(\O_U)\ar[r]^d&
\check{C}^1(\Omega^1_{U})\ar[r]^d&
\check{C}^1(\Omega^2_{U})}
\]
be the double complex and 
\[
\check{C}^0(\O_U)\os{\delta\times d}{\lra} \check{C}^1(\O_U)\times \check{C}^0(\Omega^1_{U})
\os{\cD}{\lra} \check{C}^1(\Omega^1_{U})\times
\check{C}^0(\Omega^2_{U})\os{(-d)\times\delta}{\lra}\check{C}^1(\Omega^2_{U})
\]
\[
\cD:(\alpha,\beta)\longmapsto (-d\alpha+\delta(\beta),d\beta)
\]
the associated total complex.
It gives the de Rham cohomology $H^\bullet_\dR(U)$ together with a natural map
\[
\Omega^1_S\ot H^1_\dR(U/S)\lra H^2_\dR(U),\quad dt\ot [(f)\times(z_0,z_\infty)]
\mapsto(-fdt)\times(dt\wedge z_0,dt \wedge z_\infty ).
\]
Now
\begin{align*}
\cD:\widehat{\canh}&\longmapsto (\widehat{\frac{dx}{y}}-\frac{dx}{y})
\times (d\left(\widehat{\frac{dx}{y}}\right),\frac{dxdy}{y^2})\\
&=\left(-\frac{(6g_2x^2-9g_3x-g_2^2)(g'_2x+g'_3)}{y\Delta }dt\right)\times
(d\left(\widehat{\frac{dx}{y}}\right),\frac{dxdy}{y^2})\\
&=\left(F dt\right)\times
(G_1 \frac{dtdx}{y},G_2 \frac{dtdx}{y})
\in 
\check{C}^1(\Omega^1_{U})\times\check{C}^0(\Omega^2_{U})
\end{align*}
where
\[
F:=\frac{-(\EEE)(6x^2-g_2)-\Delta'x-9g_2g_2^\prime y^2}{6y\Delta }
\]
\[
G_1:=\frac{18g_2g_2^\prime x^2+(\EEE)x
-2g_2^2g_2^\prime+9g_3g_3^\prime}{2\Delta},\quad
G_2:=\frac{g_2^\prime x+g_3^\prime}{2y^2}.
\]
This means
\[
\nabla\left(\canh\right)=dt\ot\left((-F)\times
(G_1 \frac{dx}{y},G_2 \frac{dx}{y})\right)\in \Omega^1_S\ot H^1_\dR(U/S).
\]
By Lemma \ref{A0} we get
\[
dt\ot\left((-F)\times
(G_1 \frac{dx}{y},G_2 \frac{dx}{y})\right)
\equiv dt\ot\left((-\bar{F})\times
(\bar{G}_1 \frac{dx}{y},\bar{G}_2 \frac{dx}{y})\right)
\mod \Image\check{C}^0(\O_U)
\]where
\[
\bar{F}=-\frac{\EEE}{\Delta}\frac{x^2}{y}
\]
\[
\bar{G}_1=-\frac{\Delta'}{12\Delta}+
\frac{\EEE}{\Delta}\frac{x}{2},\quad
\bar{G}_2=-\frac{\Delta'}{12\Delta}+
\frac{\EEE}{\Delta}\frac{2g_2x^2+3g_3x}{2y}
\]
and the RHS is equal to
\[
\frac{\EEE}{\Delta}dt\ot\can
-\frac{\Delta^\prime}{12\Delta}dt\ot\canh.
\]
This completes the proof of \eqref{A1-3}.

\medskip

Next we show \eqref{A1-4}.
The proof goes in the same way as above.
Let
\[
\widehat{\can}:=
(\frac{x^2}{y})\times (\frac{x}{2}\widehat{\frac{dx}{y}},\frac{(2g_2x^2+3g_3x)dx}{2y^3})
\in \check{C}^1(\O_U)\times \check{C}^0(\Omega^1_{U})
\]
be a lifting of $\can$.
Then
\begin{align*}
\cD(\widehat{\can})&=\left(
-d(\frac{x^2}{y})+\frac{x}{2}\widehat{\frac{dx}{y}}-\frac{(2g_2x^2+3g_3x)dx}{2y^3}
\right)\times
\left(d\left(\frac{x}{2}\widehat{\frac{dx}{y}}\right),
d\left(\frac{(2g_2x^2+3g_3x)dx}{2y^3}\right)\right)\\
&=(F dt)\times(G_1\frac{dtdx}{y},G_2\frac{dtdx}{y})\in 
\check{C}^1(\Omega^1_{U})
\times \check{C}^0(\Omega^2_{U})
\end{align*}
where
\[
h_1=-(g_2x+3g_3)(\EE),\quad
h_2=-4x^2(g'_2x+g'_3)\Delta
\]
\[
F=\frac{1}{8\Delta}\left(
-(6g_2g'_2x-9g_3g'_2+6g_2g'_3)y+\frac{h_1}{y}+\frac{h_2}{y^3}
\right)
-\frac{\Delta'}{12\Delta}\frac{x^2}{y}
\]
\[
G_1=\frac{30g_2g'_2x^3-(4g_2^2g'_2+9g_3g'_3)x+9(\EE)x^2-2g_2^2g'_3}{4\Delta}
\]
\[
G_2=\frac{6g_2g'_2x^3+4g'_2x^2y^2+(9g_3g'_2+6g_2g'_3)x^2+6g'_3xy^2+9g_3g'_3x}
{4y^4}
\]
This means
\[
\nabla\left(\can\right)=dt\ot\left((-F)\times
(G_1 \frac{dx}{y},G_2 \frac{dx}{y})\right)\in \Omega^1_S\ot H^1_\dR(U/S).
\]
By Lemma \ref{A0} again, we get
\begin{align*}
dt\ot\left((-F)\times
(G_1 \frac{dx}{y},G_2 \frac{dx}{y})\right)
\equiv&dt\ot\left((-\bar{F})\times
(\bar{G}_1 \frac{dx}{y},\bar{G}_2 \frac{dx}{y})\right)
\mod \Image\check{C}^0(\O_U)\\
=&\left(\frac{\Delta^\prime}{12\Delta}dt\ot\can
-\frac{g_2(\EE)}{16\Delta}dt\ot\canh\right)
\end{align*}
where
\[
\bar{F}=-\frac{\Delta'}{12\Delta}\frac{x^2}{y},\quad
\bar{G}_1=\frac{\Delta'}{12\Delta}\frac{x}{2}-
\frac{g_2(\EE)}{16\Delta}
\]
\[
\bar{G}_2=\frac{\Delta'}{12\Delta}
\frac{2g_2x^2+3g_3x}
{2y^2}
-\frac{g_2(\EE)}{16\Delta}
\]
This completes the proof of \eqref{A1-4}. QED.

\bigskip

The above computation shows the following. 
\begin{prop}\label{extraA1}
Let 
\[
\left(\frac{dx}{y}\right)_0:=\widehat{\frac{dx}{y}}+\frac{3}{2}g_2g'_2y\frac{dt}{\Delta}
=\frac{dx}{y}-\frac{(6x^2-g_2)(\EEE)+\Delta' x}{6y}\frac{dt}{\Delta}
\]
\[
\left(\frac{dx}{y}\right)_\infty:=
\frac{dx}{y}-
\frac{-g_2(\EEE)x+\Delta' x}{6y}\frac{dt}{\Delta}
\]
and 
\[
\wt{\canh}:=(0)\times\left(\left(\frac{dx}{y}\right)_0,
\left(\frac{dx}{y}\right)_\infty\right)
\in \check{C}^1(\O_U)\times \check{C}^0(\Omega^1_{U}).
\]
Let
\[
\wt{\can}=(\frac{x^2}{y})\times(\frac{x}{2}\frac{dx}{y}+F_1dt,
\frac{(2g_2x^2+3g_3x)dx}{2y^3}+F_2dt)
\in \check{C}^1(\O_U)\times \check{C}^0(\Omega^1_{U}).
\]
where
\[
F_1=-\frac{\Delta'}{12\Delta}\frac{x^2}{y}-\frac{(\EE)(g_2x+3g_3)}{8\Delta y},
\]
\[
F_2=-\frac{x^2(g'_2x+g'_3)}{2y^3}-\frac{(\EE)(g_2x+3g_3)}{8\Delta y}.
\]
Let $\cD:
\check{C}^1(\O_U)\times \check{C}^0(\Omega^1_{U})\to
\check{C}^1(\Omega^1_U)\times \check{C}^0(\Omega^2_{U})$ be as before.
Then
\begin{equation}\label{extraA2}
\cD\wt{\canh}=
\frac{\EEE}{\Delta}(dt\ot\can)'-\frac{\Delta'}{12\Delta}
(dt\ot\canh)'
\end{equation}
\begin{equation}\label{extraA21}
\cD\wt{\can}=
\frac{\Delta'}{12\Delta}(dt\ot\can)'-\frac{g_2(\EE)}{16\Delta}
(dt\ot\canh)'
\end{equation}
\begin{equation}\label{extraA3}
\cD\left(f(t)\left[\wt{\can}-\frac{\Delta'}{36(\EE)}\wt{\canh}\right]
-\frac{f'(t)\Delta}{\EEE}\wt{\canh}\right)
=(f^{\prime\prime}A+f'A'+fB)(dt\ot\omega)'.
\end{equation}
Here $A,B$ are as in Cor. \ref{PF} and we denote
\[
(dt\ot\can)'=(-\frac{x^2}{y}dt)\times(\frac{xdtdx}{2y},\frac{(2g_2x^2+3g_3x)dtdx}{2y^3})
\]
\[
(dt\ot\canh)'=(0)\times(\frac{dtdx}{y},\frac{dtdx}{y})
\]
\end{prop}


\bigskip

\noindent
Department of Mathematics, Hokkaido University,
Sapporo 060-0810,
JAPAN

\medskip

\noindent
asakura@math.sci.hokudai.ac.jp


\begin{thebibliography}{AAAi}
\bibitem[AS]{sato}
Asakura, M. and Sato, K.:
{\it Chern class and Riemann-Roch theorem for cohomology without
homotopy invariance.} (preprint), arXiv:1301.5829.
\bibitem[CDKL]{lewis}
Chen, X., Doran, C., Kerr, M. and Lewis, J.:
{\it Normal functions, Picard-Fuchs equations and elliptic fibrations on K3 surfaces.}
(preprint).
\bibitem[G]{gillet}
Gillet, H.:
{\it Riemann-Roch theorem for higher Algebraic $K$-theory.}
Adv. Math. 40 (1981), 203--289.

\bibitem[GL]{lewisJAG}
Gordon, B. and Lewis, J.: {\it Indecomposable higher Chow cycles 
on products of elliptic curves.} J. Algebraic Geom. 8 (1999), 
no. 3, 543--567. 
\bibitem[H]{h}
Hartshorne, R.: {\it On the De Rham cohomology of algebraic varieties. }
I.H.E.S. Publ. Math. No. 45 (1975), 5--99. 
\bibitem[R]{R}
Ramakrishnan, D.:
{\it Arithmetic of Hilbert-Blumenthal surfaces. }
In Number theory (Montreal, Que., 1985), 285--370, CMS Conf. Proc., 7, Amer. Math. Soc., 1987. 
\bibitem[S]{schneider}
P. Schneider:
{\it Introduction to the Beilinson conjectures.}
In Beilinson's Conjectures on Special Values of $L$-Functions
(M. Rapoport, N. Schappacher and P. Schneider, ed),
Perspectives in Math. Vol.4, 1--35, 1988.
\bibitem[Sch]{scholl}
Scholl, A.: {\it Integral elements in $K$-theory and products of modular curves.}
In The arithmetic and geometry of algebraic cycles (Banff, AB, 1998), 467--489, 
NATO Sci. Ser. C Math. Phys. Sci., 548,  (2000). 


\bibitem[Si]{Si}
   Silverman, J.: {\it Advanced topics in the arithmetic of elliptic curves}. 
   Grad.\ Texts in Math.\ 15, New York, Springer 1994.
\bibitem[AEC]{AEC}
   Silverman, J.: {\it The Arithmetic of elliptic curves}. 
   Grad.\ Texts in Math.\ 106, Springer 2009.


 
\bibitem[St]{stiller2}
   Stiller, P.:
{\it   The Picard numbers of elliptic surfaces with many symmetries.}
   Pacific J. Math.\ {\bf 128}, 157--189  (1987).
\bibitem[SZ]{SZ}
Steenbrink, J. and Zucker, S.:   
{\it Variation of mixed Hodge structure. I.} Invent. Math. 80 (1985), no. 3, 489--542.
\end{thebibliography}
\end{document}